\documentclass{article}
\usepackage{array}
\usepackage{amsmath}
\usepackage{amssymb}
\usepackage{amsthm}
\usepackage{amsfonts}
\usepackage{graphicx}
\usepackage[title]{appendix}
\usepackage[margin=0pt,font+=small,labelformat=parens,labelsep=space, skip=6pt,list=false,hypcap=false] {caption}
\usepackage{subfig}
\usepackage{amsthm}
\usepackage[top=.9in,bottom=.9in]{geometry}
 \numberwithin{equation}{section}
 
\begin{document}
%
% doesn't seem to be recognized.  \numberwithin{equation}{section}
% nor this \addtoreset{equation}{section}
%\pagestyle{empty}
% activate the following statement to obtain visible labels for each equation
%\let\labeldefs\iftrue
\let\labeldefs\iffalse
\pagestyle{plain}
%\begin{maplegroup}
%\pagenumbering{arabic}
%\newcommand{\ignore}[1]{}
\numberwithin{equation}{section}
%!!! why doesn't "\title" work???
\begin{flushleft} \vskip 0.3 in 
%%%

\centerline{Some Additions to a Family of Sums and Integrals related to Hurwitz' Zeta Function(s),}
\centerline {Euler polynomials and Euler Numbers} \vskip .3in 
\vskip .2in
\centerline{ Michael Milgram\footnote{mike@geometrics-unlimited.com}}
\centerline{Consulting Physicist, Geometrics Unlimited, Ltd.}
\centerline{Box 1484, Deep River, Ont. Canada. K0J 1P0}
\centerline{}
\centerline{May 7, 2020}
\centerline{}
revision history:

\begin{itemize}
\item{June 5, 2020.\\minor corrections}; 
\item{Feb. 3, 2021.\\changed title\\improved notation\\added section 3.1.2};
\end{itemize}

%\centerline{}
%\vskip .1in

\centerline{}
\vskip .1in

MSC classes: 	11M06, 11M35, 11M99, 26A09, 30B40, 30E20, 33C20, 33C47, 33B99, 33E20
\vskip 0.1in
\centerline{\bf Abstract}\vskip .1in

Integrals involving the kernel function ${\rm sech} (\pi x)$ over a semi-infinite range are of general interest in the study of Riemann's function $\zeta(s)$ and Hurwitz' function $\zeta(s,a)$. The emphasis here is to evaluate such integrals that include monomial moments in the integrand. These are evaluated in terms of both $\zeta(s,a)$ and its alternating equivalent $\eta(s,a)$ , thereby adding some new members to a known family of related integrals. Additionally, some finite series involving both Euler polynomials and numbers are summed in closed form.  Special cases of these results are specifically evaluated and used to verify a claimed connection with the function $\zeta(2m+1)$ for $m=1$ and $m=2$.

\section{Introduction}

Integral representations of Riemann's Zeta function $\zeta(s)$ and Hurwitz' Zeta function $\zeta(s,a)$ have a long and storied history. For example, in 1895 Jensen \cite{Jensen} presented (without proof) a result equivalent to

\begin{equation}
\displaystyle \zeta \left( s \right) ={2}^{s}\int_{0}^{\infty }\!{\frac {\sinh \left( \pi\,t \right) \sin \left( s\arctan \left( 2\,t \right)  \right) }{ \left( 4\,{t}^{2}+1 \right) ^{s/2}
\mbox{}\cosh \left( \pi\,t \right) }}\,{\rm d}t, \hspace{20pt} \Re(s)>1,
\label{Jensen}
\end{equation}
that being a specialized (a=1) case of Hermite's more general result

\begin{equation}
\displaystyle \zeta \left(s,a \right) =1/(2{a}^{s})+{\frac {{a}^{-s+1}}{s-1}}+2\,\int_{0}^{\infty }\!{\frac {\sin \left( s\arctan \left( { {t}/{a}} \right)  \right) }{ \left( {a}^{2}+{t}^{2} \right) ^{s/2} \left( {{\rm e}^{2\,\pi\,t}}-1 \right) }}\,{\rm d}t\,.
\label{Hermite}
\end{equation}

More recently, Boros, Espinosa and Moll \cite{BoEsMo}, Adamchik \cite[Section 3]{Adamchik2003} and Hu and Kim \cite{HuKim}, have studied a family of integrals 

\begin{equation}
\int_{0}^{\infty} f(t)K(a,t){\rm d}t
\label{BEM}
\end{equation}
where the kernels
\begin{align*}
K(a,t)&={1}/{(\exp(2\pi at)-1)},\\
&={1}/{(\exp(2\pi at)+1)},\\
\mbox{\rm and}\\
&={1}/{\sinh(2\pi at)}\\
\end{align*}

generalize the older results. Following that, Patkowski \cite{Patkowski2}, attempted to extend Boros, Espinosa and Moll's results to include the kernel
\begin{equation}
K(a,t)=1/\cosh(\pi t)\,.
\label{Patk}
\end{equation}
Unfortunately, that paper contains a large number of misprints (see the Appendix for corrections) as well as a fundamental error in the analysis of the purported result that incorporates kernel \eqref{Patk}. However, Patkowski's paper makes a significant contribution by showing the importance of the associated Laguerre polynomials \cite{Laguerre}

% from file Patk.mw
\iffalse 
\begin{align} \nonumber
\displaystyle {\it L} \left( n,a,x \right) &\equiv {\frac {\Gamma \left( a+1+n \right) {\mbox{$_1$F$_1$}(-n;\,a+1;\,x)}}{\Gamma \left( n+1 \right) 
\mbox{}\Gamma \left( a+1 \right) }}\\
&=\sum _{j=0}^{n}{\frac { \left( -1 \right) ^{j}\Gamma \left( n+a+1 \right) 
\mbox{}{x}^{j}}{\Gamma \left( n-j+1 \right) \Gamma \left( a+j+1 \right) \Gamma \left( j+1 \right) }}
\label{LagDef}
\end{align}
\fi

\begin{align} 
\displaystyle {\it L} \left( n,a,x \right) 
=\sum _{j=0}^{n}{\frac { \left( -1 \right) ^{j}\Gamma \left( n+a+1 \right) 
\mbox{}{x}^{j}}{\Gamma \left( n-j+1 \right) \Gamma \left( a+j+1 \right) \Gamma \left( j+1 \right) }}
\label{LagDef}
\end{align}

that arise from the identity \cite[Eq. (2.6)]{Patkowski2}, \cite[Eq. 3.769(4)]{G&R} 

\begin{equation}
\displaystyle \int_{0}^{\infty }\!{\frac {{t}^{2\,n}\sin \left( wt \right)\sin \left( s\arctan \left( {\frac {t}{a}} \right)  \right) }{ \left( {a}^{2}+{t}^{2} \right) ^{s/2}} }\,{\rm d}t=\frac { \left( -1 \right) ^{n}
\mbox{}\pi\,\Gamma \left( 2\,n+1 \right) {\it L} \left( 2\,n,s-2\,n-1,aw \right){{\rm e}^{-aw}} }{2\,\Gamma \left( s \right) 
\mbox{}{w}^{2\,n+1-s}}\,.
\label{Laguerre}
\end{equation}

Approaching from a different direction, Milgram \cite[Eq. 2.51]{MilgramEiSum} has noted an integral representation related to $\zeta(2m+1)$ that involves the kernel \eqref{Patk}, specifically

% from file EiSum.tex, taken from file Ts1.maple%

\begin{align} 
\displaystyle  \left( {2}^{2\,m+1}-1 \right) \zeta \left( 2\,m+1 \right) =A(m)
 - \left( -1 \right) ^{m}{\pi}^{2\,m}\sum _{j=0}^{m}{\frac {\mathcal{E} \left( 2\,j \right) \psi \left( 1-2\,j+2\,m \right) }{\Gamma \left( 1-2\,j+2\,m \right) \Gamma \left( 2\,j+1 \right) 
\mbox{}}}\,,
\label{Cp5a}
\end{align}
where 
% from file case_m=1.mw 

\begin{equation}
\displaystyle A \left( m \right) ={\frac {{2}^{3+2\,m}
\mbox{}{\pi}^{2\,m}}{\Gamma \left( 3+2\,m \right) }\int_{0}^{\infty }\!{\frac {{t}^{2\,m+2}{\mbox{$_3$F$_2$}(1,1,3/2;\,2+m,3/2+m;\,-4\,{t}^{2})}}{\cosh \left( \pi\,t \right) }}\,{\rm d}t}\,
\label{A(m)}
\end{equation} 
and $\mathcal{E}(2j)$ are Euler numbers (see below for a summary of notation).\newline
% see file Note_feb_7_2020.mw

It will be shown in Section \ref{sec:Eta} that the hypergeometric function in \eqref{A(m)} reduces to combinations of functions that include $\arctan(t/a)$ and $\log(1+4t^2)$; hence the interest in integrals with kernel \eqref{Patk}. Finally, Patkowski (also Glasser - private communication) demonstrates the power of Parseval's identity (q.v.) that eventually leads to the main result of this work (see Theorem 1 below), that being a related form of \cite[Theorem 4]{Patkowski2} and corrected form of \cite[Eq. (3.20)]{Patkowski2} (see also \cite[Section 3.975]{G&R}). \newline

In the above, and throughout this work 
\begin{equation}
\eta(s,a)\equiv \sum_{k=0}^{\infty} (-1)^k/(k+a)^s\,
\label{ZaltDef} 
\end{equation}
is the alternating Hurwitz function, $\Re(a)>0$  and
\begin{equation}
q\equiv \frac{a}{4}+\frac{1}{4}\,.
\label{Qid}
\end{equation}

Since much of this paper follows the approach introduced by Patkowski, we note the equivalence between the symbol $S(s,a)$ used therein, and the more fundamental function $\eta(s,a)$ used here, the transformation being (see \cite[Eq. (2.2(5))]{Sriv&Choi})

% see file Stieltjes.mw
\begin{equation}
\displaystyle S \left( s,a \right) \equiv\zeta \left(s,a \right) -\zeta \left( s,a+1/2 \right)=2^s\eta(s,2a)\,.
\label{SaDef}
\end{equation}

Throughout, the symbols $j,k,m,n,p$ are non-negative integers and $J$ is an odd integer; other symbols are complex. $\mathcal{E}(n,x)$ is an Euler polynomial, $\mathcal{E}(n)$ is an Euler number, $B(n,x)$ is a Bernoulli polynomial and  $B(n)$ is a Bernoulli number. $\gamma$ is the Euler-Mascheroni constant and $\psi(x)$ is the digamma function. $A$ and $G$ are Glaisher's and Catalan's constants respectively. Unless otherwise specified, the prime symbol $^{\prime}$ applied to $\zeta(s,a)$, $\eta(s,a)$ or $S(s,a)$ always refers to differentiation with respect to the first argument. Thus

\begin{equation}
\eta^{\,\prime}(-n,a)\equiv \frac{\partial}{\partial s}\eta(s,a)_{|s=-n}\,.
\end{equation}

For clarity, the following symbols are used interchangeably throughout

\begin{equation}
\displaystyle S^{\prime} \left( s,a \right) ={2}^{s}\ln  \left( 2 \right) \eta \left( s,2\,a \right) +{2}^{s} \eta^{\prime}   \left( s,2\,a \right) =\zeta^{\prime} \left(s,a \right) 
\mbox{}-\zeta^{\prime} \left(s,a+1/2 \right)
\label{clarity}
\end{equation}
the choice being governed by whichever variation yields typographical brevity. Also, throughout, I employ, without further comment, the identities
\begin{equation}
\displaystyle \sin \left( \arctan \left( x \right)  \right) ={\frac {x}{ \sqrt{{x}^{2}+1}}}
\label{SinId}
\end{equation}
and
\begin{equation}
\displaystyle \cos \left( \arctan \left( x \right)  \right) = \frac{1}{  \sqrt{{x}^{2}+1} }\,.
\label{CodId}
\end{equation}

% from file corollary4.mw
{\bf Theorem 1}
\iffalse 
\begin{align} \nonumber
\displaystyle \int_{0}^{\infty }\!&{\frac {{v}^{2\,n+1}\sin \left( s\arctan \left( {\frac {2v}{a}} \right)  \right)}{ \left( {a}^{2}+4\,{v}^{2} \right) ^{s/2}\cosh \left( \pi\,v \right) 
\mbox{}} }\,{\rm d}v=  {2}^{2\,n-2\,s-1} \left( -1 \right) ^{1+n}\,a\\
&\times\left( \Gamma \left( 2+2\,n \right) \sum _{j=0}^{2\,n}{\frac { \left( -a/4 \right) ^{2\,n-1-j}S \left( s-1-j,q \right) }{\Gamma \left( j+2 \right) \Gamma \left( 2\,n-j+1 \right) }}
\mbox{}+S \left( s,q \right)  \left( -a/4 \right) ^{2\,n} \right)
\label{Theorem4}
\end{align}
\fi

% see file ...\2016\New_Master_Functions_Redo_of_old_paper_from_2014\check_new_paper

\begin{equation}
\displaystyle  \int_{0}^{\infty }\!{\frac {{v}^{2\,n+1}\sin \left( s\arctan \left( {\frac {2v}{a}} \right)  \right)}{ \left( {a}^{2}+4\,{v}^{2} \right) ^{s/2}\cosh \left( \pi\,v \right) 
\mbox{}} }\,{\rm d}v=\frac{ \left( -1 \right) ^{n+1} \left( a/2 \right) ^{2\,n+1}\Gamma \left( 2+2\,n \right)}{{2}^{s}} \sum _{j=0}^{2\,n+1}\frac {\left( -2\,/{a} \right) ^{j}\,\eta \left( s-j,2\,q \right) }{\Gamma \left( 2\,n-j+2 \right) \Gamma \left( j+1 \right) } 
\label{Theorem4}
\end{equation}

\section{Proof of Theorem 1 }

Following Patkowski, the Fourier sine transform of the function $1/\cosh(\pi t/2)$ is required in order to apply Parseval's identity to evaluate the integral \cite[Theorem 4]{Patkowski2}, which transforms to \eqref{Theorem4} with even moments. Although such an identity exists \cite[Eq. 3.981.9]{G&R}, it is rather complicated, and related discussion will be deferred to Section \ref{sec:Even}.  Alternatively, consider the well-known \cite[Eq. 3.524.13] {G&R} Fourier sine transform of the associated function $t/\cosh(\pi t/2)$  

\begin{equation}
\displaystyle \int_{0}^{\infty }\!{\frac {v\sin \left( wv \right) }{\cosh \left( \pi\,v/2 \right) }}\,{\rm d}v={\frac {\sinh \left( w \right) }{  \cosh^2 \left( w \right) }}\,,
\label{sechT}
\end{equation}
from which, employing Parseval's identity along with \eqref{Laguerre} and \eqref{sechT} and choosing $s>2n-1$, we have

\begin{align} \nonumber
&\displaystyle \int_{0}^{\infty }\!\frac {{v}^{2\,n+1}\,\sin \left( s\arctan \left( {\frac {v}{a}} \right)  \right)}{ \left( {a}^{2}+{v}^{2} \right) ^{s/2}\cosh \left( \pi\,v/2 \right) 
} \,{\rm d}v={\frac { \left( -1 \right) ^{n}\Gamma \left( 2\,n+1 \right) }{\Gamma \left( s \right) }\int_{0}^{\infty }\!{\frac {{\it L} \left( 2\,n,s-2\,n-1,aw \right) \sinh \left( w \right) 
\mbox{}{{\rm e}^{-aw}}}{{w}^{2\,n+1-s}\cosh^2 \left( w \right)  }}\,{\rm d}w}\\
&=\displaystyle  \left( -1 \right) ^{n}\Gamma \left( 2\,n+1 \right) \sum _{j=0}^{2\,n}  {\frac { \left( -a \right) ^{j}
\mbox{}}{\Gamma \left( 2\,n-j+1 \right) \Gamma \left( s-2\,n+j \right) \Gamma \left( j+1 \right) }\int_{0}^{\infty }\!{\frac {{w}^{s-2\,n-1+j}\sinh \left( w \right) {{\rm e}^{-aw}}}{ \cosh^2 \left( w \right)  }}\,{\rm d}w}\, 
\label{Thm4a}
\end{align}

after applying \eqref{LagDef} and inverting the order of operations because the integral converges under the condition specified and the sum is finite. To evaluate the integral in the second equality of \eqref{Thm4a}, consider the well-known \cite[Eq. 25.11.31]{NIST} representation (and corrected form of \cite[Eq. (3.3)]{Patkowski2}), which generalizes \cite[Eq. 3.552.3]{G&R}
\begin{equation}
\displaystyle \int_{0}^{\infty }\!{\frac {{v}^{s-1}{{\rm e}^{-av}}}{\cosh \left( v \right) }}\,{\rm d}v={2}^{1-2\,s} \Gamma \left( s \right) S(s,q)\,.
\label{Eq3p3}
\end{equation}

After integration by parts, \eqref{Eq3p3} becomes

% from file Patk.mw

\begin{equation}
\displaystyle a\int_{0}^{\infty }\!{\frac {{{\rm e}^{-av}}{v}^{s}}{\cosh \left( v \right) }}\,{\rm d}v+\int_{0}^{\infty }\!{\frac {{v}^{s}{{\rm e}^{-av}}\sinh \left( v \right) }{  \cosh^2 \left( v \right)  }}\,{\rm d}v=\Gamma \left( s+1 \right) 
\mbox{}{2}^{1-2\,s} S(s,q)\,, 
\label{Thm4b}
\end{equation}

from which, using \eqref{Eq3p3} with $s:=s+1$ it follows that

\begin{equation}
\displaystyle \int_{0}^{\infty }\!{\frac {{{\rm e}^{-aw}}{w}^{s}\sinh \left( w \right) }{  \cosh^2 \left( w \right)   }}\,{\rm d}w={2}^{-2\,s}\Gamma \left( s+1 \right) 
 \left(2S(s,q)-\frac{a}{2}\,S(s+1,q) 
\mbox{} \right)\,.
\label{Thm4c}
\end{equation}
Set $s\rightarrow s-2n-1+j$ and substitute \eqref{Thm4c} into \eqref{Thm4a}, to arrive at the  following

\begin{align} \nonumber
\displaystyle \int_{0}^{\infty }&\!\frac {{v}^{2\,n+1}\sin \left( s\arctan \left( {\frac {v}{a}} \right)  \right) }{ \left( {a}^{2}+{v}^{2} \right) ^{s/2}\cosh \left( \pi\,v/2 \right) 
}\,{\rm d}v\\
&= \left( -1 \right) ^{n}{2}^{4\,n-2\,s+3}\sum _{j=0}^{2\,n}\binom{2\,n}{j} \left( -a/4 \right) ^{j} \left( S \left( s-2\,n+j-1,q \right) -\frac{a}{4}\,S \left(s-2\,n+j,q \right)  \right)\,. 
\label{Theorem4a}
\end{align} 
Since \eqref{Theorem4a} involves two separate finite sums, shift the index of one sum by unity so that the arguments of $S(s-2n+j,q)$ and $S(s-2n+j-1,q)$ coincide, and after some minor simplification, including reversal of the series, addition and subtraction of a single term, a minor change of variables in the integration term, and the use of \eqref{SaDef}, the result \eqref{Theorem4} will be found. {\bf QED}

\section {Corollaries flowing from \eqref{Theorem4}}
\subsection{Sums}
\subsubsection{A sum involving an Euler polynomial} \label{sec:EulerP}

% These can all be found in file Bernsum_3.mw%

\begin{equation}
\displaystyle \sum _{j=0}^{2\,n}{\frac { \left( -a/2 \right) ^{-j}\mathcal{E} \left( j+1,2\,q \right) }{\Gamma \left( j+1 \right) \Gamma \left( 2\,n-j+1 \right) 
\mbox{}}}={\frac {\mathcal{E} \left( 2\,n \right) }{2\,\Gamma \left( 2\,n+1 \right) {a}^{2\,n-1}}}\,.
\label{Cor1}
\end{equation}

{\bf Proof:} Set s=0 in \eqref{Theorem4a} to find
\begin{equation}
\displaystyle \sum _{j=0}^{2\,n}{\frac { \left( -a/4 \right) ^{j} \left( S \left( j-2\,n-1,q \right) -\frac{a}{4}S \left( j -2\,n,q\right)  \right) }{\Gamma \left( 2\,n-j+1 \right) 
\mbox{}\Gamma \left( j+1 \right) }}=0\,,
\label{SaSum}
\end{equation}

and invoke \cite[Lemma 3.1]{BoEsMo}, that is
\begin{equation}
\zeta(1-m,q)=-B(m,q)/m\hspace{1cm} m>0\,.
\label{Zeta2Bern}
\end{equation}

Then reverse the sum, so that, with \eqref{SaDef}, \eqref{SaSum} can be rewritten in terms of Bernoulli polynomials as
% from file Bern_Sum_3.mw

\begin{align} \nonumber
\displaystyle \sum _{j=0}^{2\,n}\frac { \left( -a/4 \right) ^{2\,n-j}}{\Gamma \left( 2\,n-j+1 \right) \Gamma \left( j+1 \right) }&\left( {\frac {B \left( j+2,q \right) }{j+2}}-{\frac {B \left( j+2,q+1/2 \right) }{j+2}} \right.\\  
&\left.
-\frac{a}{4} \left( {\frac {B \left( j+1,q \right) }{j+1}}-{\frac {B \left( j+1,q+1/2 \right) }{j+1}} \right) 
\mbox{} \right) =0\,.
\label{BernSum}
\end{align}

The third and fourth sums in \eqref{BernSum} can each be summed using \cite[Eq. 24.4.12]{NIST}, which, in terms of the variables used here reads
\begin{equation}
\displaystyle \sum _{j=0}^{2\,n}{\frac { \left( -a/4 \right) ^{-j}B \left( j+1,q \right) }{\Gamma \left( j+2 \right) \Gamma \left( 2\,n-j+1 \right) }}=-{\frac {-B \left( 2\,n+1,q-a/4 \right) 
\mbox{}+ \left( -a/4 \right) ^{2\,n+1}}{\Gamma \left( 2+2\,n \right)  \left( -a/4 \right) ^{2\,n}}}\,.
\label{BernId1}
\end{equation}

The first two sums in \eqref{BernSum} can be reduced by reference to \cite[Eq. 24.4.23]{NIST}, that being
\begin{equation}
\displaystyle {B} \left( j+2,q+1/2 \right) -{B} \left( j+2,q \right) ={\frac { \left( j+2 \right) \mathcal{E} \left( j+1,2\,q \right) 
\mbox{}}{{2}^{j+2}}}\,.
\label{BernId2}
\end{equation}
In the case that $q=1/4$ %along with \cite[Eq. 24.4.31]{NIST} which in turn reads
\begin{equation}
\displaystyle B \left( 2\,n+1,1/4 \right) =-B \left( 2\,n+1,3/4 \right) =-{\frac { \left( 2\,n+1 \right) {\mathcal{E}} \left( 2\,n \right) }{{4
\mbox{}}^{2\,n+1}}}
\label{BernId3}
\end{equation}
which, applied to \eqref{BernId1} along with the identity

\begin{equation}
\displaystyle \mathcal{E} \left( n,1/2 \right) ={2}^{-n}\mathcal{E} \left( n \right)
\label{EuId}
\end{equation}

eventually yields \eqref{Cor1}. {\bf Q.E.D.}  
\newline

{\bf Remarks:}  
\begin{itemize}
\item

Note that \eqref{Cor1} differs from the well-known result \cite[Eq. (24.4.13)]{NIST} which, in the present notation, reads

% see file Corollary4.mw for some proofs of this; also file Bern_Sum_3.mw

\begin{equation}
\displaystyle \sum _{j=0}^{2\,n}{\frac {\left( -a/2 \right) ^{-j}\,\mathcal{E} \left( j,2\,q \right)  }{\Gamma \left( j+1 \right) \Gamma \left( 2\,n-j+1 \right) }}=
\mbox{}{\frac {\mathcal{E} \left( 2\,n \right) }{\Gamma \left( 1+2\,n \right) {a}^{2\,n}}}\,;
\label{Euid1}
\end{equation}
\item
\eqref{Euid1} can be obtained from \eqref{Theorem4} by setting $s=0$ and employing the above identities. See \cite[Section 5] {NaRuVa} for similar results;
\item
The result \eqref{SaSum} is equivalent to
\begin{equation}
\displaystyle \sum _{j=0}^{2\,n+1}\left( \frac{-2}{a} \right) ^{j}{\frac {\eta \left( -j,2\,q \right) }{\Gamma \left( 2\,n-j+2 \right) \Gamma \left( j+1 \right) } }=0\,.
\label{Sx}
\end{equation}
\end{itemize}

\subsubsection{A sum involving Euler Numbers} \label{sec:EulerN}
If $J$ is an odd integer,
\begin{align} \nonumber
\displaystyle &\sum _{j=0}^{1+n}\left( \frac{-1}{2J}\right)^{j}\frac {\mathcal{E} \left( 1+n+j \right) }{\Gamma \left( n-j+2 \right) \Gamma \left( j+1 \right) 
\mbox{}} ={\frac {2\, \left( -1 \right) ^{J/2+1/2}}{ \left( 2\,J \right) ^{1+n}\Gamma \left( 2+n \right) }\sum _{k=1}^{J/2-1/2}{\frac { \left( -1 \right) ^{k}}{ \left( {J}^{2}-4\,{k}^{2} \right) ^{-1-n}}}}\\
&-{2}^{2+n}\sum _{j=0}^{1+n} \left( -\frac{1}{J} \right) ^{j} \frac {1}{\Gamma \left( n-j+2 \right) \Gamma \left( j+1 \right) }\sum _{k=1}^{J} \left( -1 \right) ^{k} \left( k-1/2 \right) ^{1+n+j} +{\frac { \left( -1 \right) ^{J/2+1/2} \left( J/2 \right) ^{1+n}}{\Gamma \left( 2+n \right) }}.
\label{New3c}
\end{align}

{\bf Proof:}
% from file 2016\\new_Masterfunction\\redo_of_old_paper_from 2014\\Example_2c_continued.mw

From Glasser and Milgram, \cite[Eq. (3.9)]{Master}, let $2k-1=J$, and after an obvious change of integration variables, obtain the identity

\begin{align} \nonumber
\displaystyle \int_{0}^{\infty }&\!{\frac {{v}^{-s} \left( \cos \left( s\pi/2 \right) \cos \left( s\arctan \left( {\frac {v}{J}} \right)  \right) +\sin \left( s\pi/2 \right) \sin \left( s\arctan \left( {\frac {v}{J}} \right)  \right) 
\mbox{} \right)}{\left( {J}^{2}+{v}^{2} \right) ^{s/2}\cosh \left( \pi\,v \right) }  }\,{\rm d}v\\
&={2}^{2\,s} \left( -1 \right) ^{3/2+J/2}\sum _{k=1}^{J/2-1/2} \frac{\left( -1 \right) ^{k}}{   \left( {J}^{2}-4\,{k}^{2} \right) ^{s}}  
\mbox{}-\frac{ \left( -1 \right) ^{J/2+1/2}}{2} \left(\frac{2}{J} \right) ^{2\,s}\,.
\label{Rkb}
\end{align}

In \eqref{Rkb}, let $s=-2n-1$ to give

\begin{align}  \nonumber
\displaystyle \int_{0}^{\infty }\!{\frac {{v}^{1+2\,n}\,\sin \left(  \left( 1+2\,n \right) \arctan \left( {\frac {v}{J}} \right)  \right)}{ \left( {J}^{2}+{v}^{2} \right) ^{-n-1/2}
\mbox{}\cosh \left( \pi\,v \right) } }\,{\rm d}v&={\frac { \left( -1 \right) ^{3/2+J/2+n}}{{2}^{4\,n+2}}\sum _{k=1}^{J/2-1/2}{\frac { \left( -1 \right) ^{k}}{ \left( {J}^{2}-4\,{k}^{2} \right) ^{-1-2\,n}}}}\\
&+{\frac { \left( -1 \right) ^{3/2+J/2+n}{J}^{4\,n+2}}{{2}^{4\,n+3}}}
\label{Rkbo}
\end{align}
and, after setting $s=-1-2n$ and $a=2J$ in \eqref{Theorem4} which now becomes 

\begin{align} \nonumber
\displaystyle \int_{0}^{\infty }\!\frac {{v}^{1+2\,n}\,\sin \left(  \left( 1+2\,n \right) \arctan \left( {\frac {v}{J}} \right)  \right)}{ \left( {J}^{2}+{v}^{2} \right) ^{-n-1/2}
\mbox{}\cosh \left( \pi\,v \right) }& \,{\rm d}v\\
&\hspace{-70pt}= \left( -1 \right) ^{n}{J}^{1+2\,n}\,\Gamma \left( 2+2\,n \right)\sum _{j=0}^{1+2\,n}\left( \frac{-1}{J} \right) ^{j}{\frac {\eta \left( -1-2\,n-j,J+1/2 \right) }{\Gamma \left( 2\,n-j+2 \right) \Gamma \left( j+1 \right) } }
\label{T1A}
\end{align}
compare the right-hand sides of both, yielding a closed form for a particular sum of alternating Hurwitz zeta functions $\eta(s,a)$. That is

\begin{align} \nonumber
\displaystyle \sum _{j=0}^{1+2\,n} \left(\frac{-1}{J}\right) ^{j}{\frac {\eta \left( -1-2\,n-j,J+1/2 \right) }{\Gamma \left( 2\,n-j+2 \right) \Gamma \left( j+1 \right) 
\mbox{}}}&=\,-{\frac {{2}^{-2-4\,n} \left( -1 \right) ^{J/2+1/2}
\mbox{}}{{J}^{1+2\,n}\Gamma \left( 2+2\,n \right) }\sum _{k=1}^{J/2-1/2}{\frac { \left( -1 \right) ^{k}}{ \left( {J}^{2}-4\,{k}^{2} \right) ^{-1-2\,n}}}}\\
&-{\frac {{J}^{1+2\,n} \left( -1 \right) ^{J/2+1/2}{2}^{-4\,n-3}}{\Gamma \left( 2+2\,n \right) }}\,.
\label{New1A}
\end{align}

At this point, it is instructive to skip ahead to Section \ref{sec:T2andProof} where another equivalent set of calculations could be performed, this time setting $s=-2n$ in the identities corresponding to the above. Left as an exercise for the reader, an identity corresponding to \eqref{New1A} will be found that becomes identical to \eqref{New1A} after setting $n:=n+1/2$. This means that \eqref{New1A} is valid for both even and odd values of $n$, so, setting $2n:=n-1$ in \eqref{New1A} yields a more general form and an intermediate result of interest, that being

\begin{align} \nonumber
\displaystyle\sum _{j=0}^{n} \left(\frac{-1}{J} \right) ^{j}{\frac {\eta \left( -n-j\,,J+1/2 \right) }{\Gamma \left( n+1-j \right) \Gamma \left( j+1 \right) }}=\,&-{\frac { \left( -1 \right) ^{J/2+1/2}}{ \left( 4\,J \right) ^{n}\Gamma \left( 1+n \right) }\sum _{k=1}^{J/2-1/2}{\frac { \left( -1 \right) ^{k}}{ \left( {J}^{2}-4\,{k}^{2} \right) ^{-n}}}}
\mbox{}\\
&-{\frac { \left( J/4 \right) ^{n} \left( -1 \right) ^{J/2+1/2}}{2\,\Gamma \left( 1+n \right) }},
\label{New3a}
\end{align}
where, it is reiterated, $J$ is an odd integer. Since the first argument of $\eta$ in \eqref{New3a} is always a negative integer, that result can be further reduced by noting the identities \cite[Eq. (2.7)]{HuKimKim}

\begin{equation}
\eta(-m,z)=\mathcal{E}(m,z)/2
\label{EtaMz}
\end{equation}
and
\begin{equation}
\mathcal{E}(m,z)=2(z-1)^m-\mathcal{E}(m,z-1),
\label{Ezm1}
\end{equation}
the latter yielding a more general form

\begin{equation}
\displaystyle {\mathcal{E}} \left( m,z \right) =-2\,\sum _{k=1}^{K} \left( -1 \right) ^{k} \left( z-k \right) ^{m}+ \left( -1 \right) ^{K}{\mathcal{E}} \left( m,z-K \right) 
\label{EuRecur}
\end{equation}
where $K$ is a positive integer. To \eqref{New3a}, first apply \eqref{EtaMz}, then set $z=J+1/2$ and $K=J$ in \eqref{EuRecur} and, after some simplification, including the use of \eqref{EuId}, the result \eqref{New3c} will emerge. {\bf QED}.\newline

{\bf Special Case J=1}

\begin{equation}
\displaystyle \sum _{j=0}^{n}{\frac { \left( -1/2 \right) ^{j}{\mathcal{E}} \left( n+j \right) }{\Gamma \left( n+1-j \right) \Gamma \left( j+1 \right) 
\mbox{}}}={\frac {{2}^{-n}}{\Gamma \left( 1+n \right) }}\,.
\label{ScJ1}
\end{equation}

{\bf Remark:} For large values of $n$, the left hand side of \eqref{New3c} can lose many significant digits in a numerical evaluation.

\subsection{Corollary \ref{sec:3.2}}\label{sec:3.2}

% This is from file Theorem4.mw, later improved in file Corollary4.mw
\begin{align} \nonumber
\displaystyle \int_{0}^{\infty }\!&{\frac {{v}^{2+2\,n}}{ \left( {a}^{2}+4\,{v}^{2} \right) \cosh \left( \pi\,v \right) }}\,{\rm d}v= \left( -1 \right) ^{n}{2}^{2\,n-3} \left( a/4 \right) ^{2\,n}\\
&\times \left( \Gamma \left( 2+2\,n \right) \sum _{j=0}^{2\,n}{\frac { \left(- 2/{a} \right) ^{j}\mathcal{E} \left( j,2\,q \right) }{\Gamma \left( j+2 \right) \Gamma \left( 2\,n-j+1 \right)
\mbox{}}}-\frac{a}{2} \left( \psi \left( q+1/2 \right) -\psi \left( q \right)  \right) 
\mbox{} \right)  ,\hspace{.1cm} n\geq 0 \,.
\label{Cor2a}
\end{align}

{\bf Proof:}

In \eqref{Theorem4}, set $s=1$ noting the limit \cite[Eq. 25.11.31]{NIST}
\begin{equation}
\displaystyle \lim _{s\rightarrow 1}S \left( s,q \right) =\psi \left( q+1/2 \right) -\psi \left( q \right) \,.
\label{Lims1}
\end{equation}

Apply the identity \eqref{Zeta2Bern} and after some simplification, find

\begin{align} \nonumber
\displaystyle &\int_{0}^{\infty }\!{\frac {{v}^{2+2\,n}}{ \left( {a}^{2}+4\,{v}^{2} \right) \cosh \left( \pi\,v \right) }}\,{\rm d}v= {2}^{2\,n-4} \left( -1 \right) ^{n+1}\left( a/4 \right) ^{2\,n}\\
&\times \left( 4\,\Gamma \left( 2+2\,n \right)\sum _{j=0}^{2\,n}{\frac {\left( B\left( j+1,q \right) -B \left( j+1,q+1/2 \right)\right)\left( -4/a \right) ^{j} }{ \left( j+1 \right) \Gamma \left( j+2 \right) 
\mbox{}\Gamma \left( 2\,n-j+1 \right) } } +a \left( \psi \left( q+1/2 \right) -\psi \left( q \right)  \right)  \right) 
\,.
\label{Cor2b}
\end{align}

The difference of terms involving Bernoulli polynomials can be rewritten with reference to \eqref{BernId2} after which \eqref{Cor2b} is easily reduced to \eqref{Cor2a}. {\bf QED}\newline

Special cases: 
\begin{itemize}

\item
$a=1$

% This is from file Theorem4.mw, later improved in file Corollary4.mw; see also Cor4_Special_cases.mw 

\begin{align} \nonumber
\displaystyle \int_{0}^{\infty }\!&\frac {{v}^{2+2\,n}}{ \left( 4\,{v}^{2}+1 \right) \cosh \left( \pi\,v \right) }\,{\rm d}v\\
& =\left( -1 \right) ^{n}{2}^{-2-2\,n}
\left( n+(1-\ln  \left( 2 \right))/2+\Gamma \left( 2+2n \right) \sum _{j=0}^{n-1}{\frac {{2}^{2\,j}\,\mathcal{E} \left( 2\,j+1,0 \right) }{\Gamma \left( 3+2j \right) \Gamma \left( 2\,n-2\,j \right) }}
\mbox{} \right)\,; 
\label{C2na1}
\end{align}
\item$a=2$

\begin{align} \nonumber
\displaystyle \int_{0}^{\infty }\!&{\frac {{v}^{2\,n+2}}{ \left( {v}^{2}+1 \right) \cosh \left( \pi\,v \right) }}\,{\rm d}v\\
&= \left( -1 \right) ^{n+1} \left(\frac{\Gamma \left( 2\,n+2 \right)}{2} \sum _{j=1}^{n}{\frac {\mathcal{E} \left( 2\,j \right) {2}^{-2\,j}}{\Gamma \left( 2\,j+2 \right) \Gamma \left( 2\,n-2\,j+1 \right) }}
\mbox{}-\pi/2+{2}^{-2\,n}+n+1/2 \right) \,;
\label{C2na2}
\end{align}

\item $a=4$
\begin{align} \nonumber
\displaystyle \int_{0}^{\infty }\!{\frac {{v}^{2\,n+2}}{ \left( {v}^{2}+4 \right) \cosh \left( \pi\,v \right) }}\,{\rm d}v
&= \left( -1 \right) ^{n} \left( {2}^{2\,n-1}\Gamma \left( 2\,n+2 \right)\sum _{j=1}^{n}{\frac {{2}^{-4j}\,\mathcal{E} \left( 2\,j \right) }{\Gamma \left( 2\,j+2 \right) \Gamma \left( 2\,n-2\,j+1 \right) }}\right.\\
&\left. - \left( \pi-n-1/2 \right) {2}^{2\,n}
\mbox{}+{2}^{-2\,n} \left( {3}^{2\,n+1}-1/3 \right)  \right)\,.
\label{C2na4}
\end{align}

\end{itemize}
% these all come from file Cor4 Special Cases.mw

{\bf Remark:} Note that \eqref{Cor2b} extends results found in \cite[Section 3.522]{G&R}. Simplification of the above special cases employ 
\eqref{Ezm1} as well as \cite[Eq. 24.5.2]{NIST} and \cite[Eq. 24.4.13]{NIST}.

\subsection{Corollary 3.3}\label{sec:3.3}

% these all  came from file Theorem4.mw, updated in file Corollary4.mw and "Check New paper.mw"

\begin{equation}
\displaystyle\int_{0}^{\infty }\!{\frac {{v}^{2\,n+1}\,\arctan \left( {\frac {2v}{a}} \right)}{\cosh \left( \pi\,v \right) } }{\rm d}v= \left( -1 \right) ^{n+1}
 \left( \frac{a}{2} \right) ^{2\,n+1}\Gamma \left( 2+2n \right) \sum _{j=0}^{2\,n+1}\left(\frac{-4}{a}  \right) ^{j}{\frac {\zeta^{\prime} \left(-j,q \right) -\zeta^{\prime} \left(-j\,,q+1/2 \right) }{\Gamma \left( 2\,n-j+2 \right) \Gamma \left( j+1 \right) 
\mbox{}} }
\label{Tc3a}
\end{equation}

{\bf Proof:}

% These will all be found in file Corollary4.mw%

In \eqref{Theorem4}, operate with ${ \frac{\partial}{\partial s}}$ and set $s=0$, noting that the application of \eqref{Sx} results in most terms vanishing, and \eqref{Tc3a} immediately follows. It is worth noting \cite[Eq. 25.11.18]{NIST} that
\begin{equation}
\zeta^{\prime}(0,q)=\ln(\Gamma(q))-\frac{1}{2}\ln(2\pi)\,.
\label{Zp0}
\end{equation} {\bf QED}.
\newline

The following special cases can profitably be compared with examples given in Section 5 of \cite{BoEsMo}.\newline

{\bf Special case $n=0$.} \newline

{\bf Note:} None of the following cases appear to be known to Mathematica \cite{Math}.

\begin{equation}
\displaystyle \int_{0}^{\infty }\!{\frac {v\,\arctan \left( {\frac {2\,v}{a}} \right)}{\cosh \left( \pi\,v \right) } }\,{\rm d}v=a/2 \left( -\ln  \left( \Gamma \left( q \right)  \right) +\ln  \left( \Gamma \left( q+1/2 \right)  \right)  \right) 
\mbox{}+2\,(\zeta^{\prime} \left(-1,q \right) -\zeta^{\prime} \left( -1,q+1/2 \right)) 
\label{Cn0}
\end{equation}
\begin{itemize}
\item

If $a=1$ then \eqref{Cn0} reduces to

% case a=1 is done in file Patk.mw

\begin{equation}
\displaystyle \int_{0}^{\infty }\!{\frac {v\arctan \left( 2\,v \right) }{\cosh \left( \pi\,v \right) }}\,{\rm d}v=-3\,\zeta^{\prime} \left( -1 \right) +\ln  \left( 2 \right)/6
\mbox{}-\ln  \left( 2\,\pi \right)/4 
\label{a0n1}
\end{equation}
and by writing $\zeta^{\prime} \left(-1 \right)$ in the form of its underlying sum, we identify \cite{Maple} 
\begin{equation}
\displaystyle \zeta^{\prime} \left(-1 \right) =-\gamma/12-\ln  \left( 2\,\pi \right)/12 +1/12+{\frac {\zeta^{\prime} \left(2 \right) }{{2\pi}^{2}}}\,.
\label{Z1m1}
\end{equation}
Furthermore, from \cite{Math}
\begin{equation}
\displaystyle \zeta^{\prime} \left( 2 \right) ={\pi}^{2} \left( \gamma+\ln  \left( 2\,\pi \right) -12\,\ln  \left( A \right)  \right)/6\,, 
\label{Z1p2}
\end{equation}
all of which identifies
\begin{equation}
\displaystyle \int_{0}^{\infty }\!{\frac {v\arctan \left( 2\,v \right) }{\cosh \left( \pi\,v \right) }}\,{\rm d}v=-\,\ln  \left( \pi \right)/4 -1/4+3\,\ln  \left( A \right) 
\mbox{}-\ln  \left( 2 \right)/12 \,.
\label{C3n0a1}
\end{equation}
\item

% from  file Theorem4.mw near Eqs. 2.3.4

If $a=2$ then \eqref{Cn0} reduces to
\begin{equation}
\displaystyle \int_{0}^{\infty }\!{\frac {v\arctan \left( v \right) }{\cosh \left( \pi\,v \right) }}\,{\rm d}v=-2\,\ln  \left( \Gamma \left( 3/4 \right)  \right) -\ln  \left( 2 \right)/2 
\mbox{}+\ln  \left( \pi \right) +2\,\zeta^{\prime} \left( -1,3/4 \right) -2\,\zeta^{\prime} \left( -1,1/4 \right)\,.
\label{N0a2}
\end{equation}

From Miller and Adamchik \cite[Eq. (5)]{Miller&Adamchik} (also \cite[Eq. 25.11.21]{NIST})
\begin{equation}
\displaystyle \zeta^{\prime} \left(-1,3/4 \right) -\zeta^{\prime} \left( -1,1/4 \right) =-{\frac {G}{4\pi}}
\mbox{}+{\frac {\psi^{\prime} \left( 3/4 \right) }{32\,\pi}}-\pi/32
\label{Zdiff}
\end{equation}
and again from those same sources
\begin{equation}
\displaystyle \psi^{\prime} \left( 3/4 \right) =  \zeta \left(2,3/4 \right) ={\pi}^{2}-8\,G 
\label{Psi3Q}
\end{equation}
all of which leads to
\begin{equation}
\displaystyle \int_{0}^{\infty }\!{\frac {v\arctan \left( v \right) }{\cosh \left( \pi\,v \right) }}\,{\rm d}v=-\frac{1}{2}\,\ln  \left( 2 \right) -2\,\ln   \Gamma \left( 3/4 \right)  
\mbox{}+\ln  \left( \pi \right) -{ {{G}}/{\pi}}\,.
\label{C3n0a2}
\end{equation}
\item
If $a=4$ then \eqref{Cn0} reduces to

\begin{equation}
\displaystyle \int_{0}^{\infty }\!{\frac {v\arctan \left( v/2 \right) }{\cosh \left( \pi\,v \right) }}\,{\rm d}v=2\,\ln  \left( {\frac {\Gamma \left( 7/4 \right) }{\Gamma \left( 5/4 \right) }} \right) 
\mbox{}+2\,\zeta^{\prime} \left(-1,1/4 \right) -\ln  \left( 2 \right) -2\,\zeta^{\prime} \left(-1,3/4 \right) -3/2\,\ln  \left( 3/4 \right) 
\label{n0a4}
\end{equation}
where the recursion
\begin{equation}
\displaystyle \zeta^{\prime} \left(-1,q+1/2 \right) =\zeta^{\prime} \left(-1,q-1/2 \right) + \left( q-1/2 \right) \ln  \left( q-1/2 \right) 
\label{zm1}
\end{equation}
has been used to convert the original reduction of \eqref{Cn0} into a range where \cite[Eq. (5)]{Miller&Adamchik} is valid when $q>1$. Following the use of \eqref{Zdiff} and \eqref{Psi3Q} we eventually find

\begin{equation}
\displaystyle \int_{0}^{\infty }\!{\frac {v\arctan \left( v/2 \right) }{\cosh \left( \pi\,v \right) }}\,{\rm d}v=\ln  \left( 2 \right) +4\,\ln  \left( \Gamma \left( 3/4 \right)  \right) 
\mbox{}+\frac{1}{2}\ln  \left( 3 \right) -2\,\ln  \left( \pi \right) +{ {{G}/}{\pi}}\,.
\label{C3n0a4}
\end{equation}

% from file Theorem4.maple

\item
The above results can be associated with well-known arctangent identities to obtain other interesting identities. From \cite{MathWorldInvTan} we find two identities valid for all values of their arguments, those being
\begin{align}
\displaystyle &\arctan \left( v/2 \right) =\arctan \left( v \right) -\arctan \left( {\frac {v}{{v}^{2}+2}} \right) 
\label{Actv2}\\
\displaystyle &\arctan \left( v \right) =2\,\arctan \left( v/2 \right) -\arctan \left( {\frac {{v}^{3}}{3\,{v}^{2}+4}} \right) 
\label{Actv3}
\end{align}

Substitution into \eqref{C3n0a2} and \eqref{C3n0a4} yields
\begin{equation}
\displaystyle \int_{0}^{\infty }\!{\frac {\,v\,\arctan \left( {\displaystyle \frac {v}{{v}^{2}+2}} \right)}{\cosh \left( \pi\,v \right) } }\,{\rm d}v=3\,\ln  \left( \pi \right) -6\,\ln  \left( \Gamma \left( 3/4 \right)  \right) 
\mbox{}-3/2\,\ln  \left( 2 \right) -\ln  \left( 3 \right)/2 -{\frac {2\,G}{\pi}}
\label{AtanId1}
\end{equation}
and
\begin{equation}
\displaystyle \int_{0}^{\infty }\!{\frac {v\,\arctan \left( {\displaystyle \frac {{v}^{3}}{3\,{v}^{2}+4}} \right)}{\cosh \left( \pi\,v \right) } }\,{\rm d}v=10\,\ln  \left( \Gamma \left( 3/4 \right)  \right) 
\mbox{}+5\,\ln  \left( 2 \right)/2 -5\,\ln  \left( \pi \right) +{\frac {3\,G}{\pi}}+\ln  \left( 3 \right) \,.
\label{AtanId2}
\end{equation}
\end{itemize}

{\bf Special case $n=1$}\newline

\begin{align} 
\displaystyle \int_{0}^{\infty }\!{\frac {{v}^{3}\arctan \left( {\frac {2\,v}{a}} \right) }{\cosh \left( \pi\,v \right) }}\,{\rm d}v=&\frac{a^3}{8}\ln  \left( {\frac {\Gamma \left( q \right) }{\Gamma \left( q+1/2 \right) }} \right) 
\mbox{}-\frac{3a^2}{2}S^{\prime}\left(-1,q \right) 
+6\,a\,S^{\prime} \left(-2,q \right)  -8\,S^{\prime} \left(-3,q \right) \,\,.
\label{Cn1}
\end{align}
\begin{itemize}
\item
% from file Theorem4.maple

If $a=1...$, we find

\begin{align} \nonumber
\displaystyle \int_{0}^{\infty }\!{\frac {{v}^{3}\arctan \left( 2\,v \right) }{\cosh \left( \pi\,v \right) }}\,{\rm d}v&=3/2\,\left(\zeta^{\prime} \left(-1 \right)-\zeta^{\prime} \left(-1,1/2 \right)\right) 
\mbox{} +6\,\zeta^{\prime} \left( -2,1/2 \right) -6\,\zeta^{\prime} \left(-2 \right) \\
&-8\,\zeta^{\prime} \left(-3,1/2 \right) +8\,\zeta^{\prime} \left(-3 \right) 
\mbox{}+\ln  \left( \pi \right)/16. 
\label{N1a1}
\end{align}

From Mathematica \cite{Math}
% from file Theorem4.mw

\begin{align} 
\displaystyle \zeta^{\prime} \left( -1,1/2 \right) &=-\,\ln  \left( 2 \right)/24 -1/24+\ln  \left( A \right)/2 
\label{Zm1h}\\
\displaystyle \zeta^{\prime} \left(-2,1/2 \right) &={\frac {3\,\zeta \left( 3 \right) }{{16\,\pi}^{2}}}
\label{Zm2h}\\
\displaystyle \zeta^{\prime} \left(-3,1/2 \right) &={\frac {\ln  \left( 2 \right) }{960}}-{\frac {7\,\zeta^{\prime} \left( -3 \right) }{8}}
\label{Zm3h}
\end{align}
and from \cite{WikiPval},
\begin{equation}
\zeta^{\prime}(-2)=-\zeta(3)/(4\pi^2)\,.
\label{Zpm2}
\end{equation}
Together with other identifications noted previously, after simple calculation, we obtain

\begin{equation}
\displaystyle \int_{0}^{\infty }\!{\frac {{v}^{3}\arctan \left( 2\,v \right) }{\cosh \left( \pi\,v \right) }}\,{\rm d}v={\frac {13\,\ln  \left( 2 \right) }{240}}-9\,\ln  \left( A \right)/4 
\mbox{}+\ln  \left( \pi \right)/16 +15\,\zeta^{\prime} \left(-3 \right) +3/16+{\frac {21\,\zeta \left( 3 \right) }{8\,{\pi}^{2}}}
\label{N1a1f}
\end{equation} 
From \cite[Section 5]{Adamchik1998} we have
\begin{equation}
\zeta^{\prime} \left(-3 \right)=-\,\frac{11}{720}-\ln\left(A_{3}\right)\,,
\end{equation} 
where $A_{3}$ is a generalized Glaisher constant.
\newline

\item

% from file Theorem4.maple

If $a=2 $, we find
\begin{align} \nonumber
\displaystyle \int_{0}^{\infty }\!{\frac {{v}^{3}\arctan \left( v \right) }{\cosh \left( \pi\,v \right) }}\,{\rm d}v&= 6\left( \zeta^{\prime} \left( -1,1/4 \right) -\zeta^{\prime} \left(-1,3/4 \right)  \right) +8\left( \zeta^{\prime} \left(-3,1/4 \right) -\zeta^{\prime} \left(-3,3/4 \right)  \right) 
\mbox{}\\
&-12 \left(\zeta^{\prime} \left(-2,1/4 \right) -\zeta^{\prime} \left( -2,3/4 \right)  \right) 
\mbox{}-7/4\,\ln  \left( 2 \right) +\ln  \left( {\frac { 2 \sqrt{\,2}\left( \Gamma \left( 3/4 \right)  \right) ^{2}}{\pi}} \right)\,. 
\label{N1a2}
\end{align}

In addition to the identities utilized previously, we also have, from \cite[Eq. (5)]{Miller&Adamchik}, \cite[Eq. (6)]{Miller&Adamchik} and \cite[Eq. (12)]{Miller&Adamchik} respectively,
\begin{equation}
\zeta^{\prime} \left(-3,1/4 \right) -\zeta^{\prime} \left(-3,3/4 \right) ={\frac {\pi}{128}}-{\frac {\psi^{(3)} \left(1/4 \right) }{1024\,{\pi}^{3}}}\,,
\label{V4mV5}
\end{equation}

\begin{equation}
\displaystyle \zeta^{\prime} \left( -2,1/4 \right) -\zeta^{\prime} \left( -2,3/4 \right) =-\,\ln  \left( \pi \right)/32 -{\frac {3\,\ln  \left( 2 \right) }{32}}
\mbox{}+{\frac{3}{64}}-\gamma/32+{\frac {\zeta^{\prime} \left( 3,1/4 \right) -\zeta^{\prime} \left( 3,3/4 \right) }{64\,{\pi}^{3}}}\,,
\label{Rad1}
\end{equation}
and
\begin{equation}
\displaystyle \zeta^{\prime} \left(3,3/4 \right) =120\,\ln  \left( 2 \right) \zeta \left( 3 \right) +56\,\zeta^{\prime} \left( 3 \right) -\zeta^{\prime} \left( 3,1/4 \right) \,.
\label{Mult1}
\end{equation}
In the above, $\psi^{(3)} \left(1/4 \right)$, the polygamma function, can also be written 
\begin{equation}
\psi^{(3)} \left(1/4 \right)=8\pi^{4}+768\beta(4)
\label{TriPsi}
\end{equation}
where $\beta(4)$ is Dirichlet's Beta function (see \cite{Kolbig1996}). These identities eventually yield
\begin{align} \nonumber
\displaystyle \int_{0}^{\infty }\!{\frac {{v}^{3}\arctan \left( v \right) }{\cosh \left( \pi\,v \right) }}\,{\rm d}v&=2\,\ln  \left( \Gamma \left( 3/4 \right)  \right) 
\mbox{}+{ {7\,\ln  \left( 2 \right) }/{8}}-5\ln  \left( \pi \right)/8 +3\gamma/8-{{9}/{16}}+{ {3\,G}/{\pi}}
\\
&+\frac{1}{\pi^{3}} \left({{{45\,\ln  \left( 2 \right) \zeta \left( 3 \right) }/{2}}+21\,\zeta^{\prime} \left(3 \right)/2 -3\,\zeta^{\prime} \left(3,1/4 \right)/8 -6\,\beta \left( 4 \right)
\mbox{}}{}\right)\,,
\label{N1a1b}
\end{align}
\end{itemize}
which is likely the simplest form that exists in terms of independent fundamental constants (see \cite{RAMMURTY2007}, \cite{Waldschmidt} and \cite[Eq. (21)]{Miller&Adamchik}). Results equivalent to \eqref{AtanId1} and \eqref{AtanId2} can now be obtained using \eqref{Actv2} and \eqref{Actv3}.

\subsection{Corollary 3.4} \label{sec:3.4}

% from file corollary4a.mw Eq. 11

\begin{align} \nonumber
&\displaystyle a\int_{0}^{\infty }\!{\frac {{v}^{2\,n+1}\arctan \left( {\frac {2v}{a}} \right)}{ \left( {a}^{2}+4\,{v}^{2} \right) \cosh \left( \pi\,v \right) } }\,{\rm d}v-\int_{0}^{\infty }\!{\frac {{v}^{2\,n+2}\ln  \left( {a}^{2}+4\,{v}^{2} \right) }{ \left( {a}^{2}+4\,{v}^{2} \right) \cosh \left( \pi\,v \right) }}\,{\rm d}v
\mbox{}= 
 \left( -{a}^{2}/4 \right) ^{n}\left[ \frac{a\ln  \left( 2 \right) \left( \psi \left( q+1/2 \right) -\psi \left( q \right)  \right)}{8}\right.\\
&\left. -\frac{a\, \eta^{\prime}(1,2q)}{4}+\frac{\Gamma \left( 2\,n+2 \right)}{2} \left(\sum _{j=0}^{2\,n}{\frac {S^{\prime}\!   \left(-j,q \right)\left( -4/a \right) ^{j} }{\Gamma \left( j+2 \right) \Gamma \left( 2\,n-j+1 \right) } }
\mbox{}-{ \ln  \left( 2 \right)}{} \sum _{j=0}^{2\,n}{\frac { \mathcal{E} \left( j,2\,q \right)\left( -2/a \right) ^{j} }{\Gamma \left( j+2 \right) \Gamma \left( 2\,n-j+1 \right) } }
\mbox{} \right) \right]
\label{Cor4}
\end{align}
where
\begin{equation}
 \eta^{\prime}(s,2q)\equiv-\,\sum _{j=0}^{\infty }{\frac { \left( -1 \right) ^{j}\ln  \left( j+2\,q \right) }{(j+2\,q)^s}}
\label{Zalt}
\end{equation}

is the first derivative of the alternating Hurwitz zeta function with respect to $s$ (see \eqref{ZaltDef}).\newline

{\bf Proof:}\newline

% from file Stieltjes.mw

Differentiate \eqref{Theorem4} with respect to $s$, reverse the resulting sums and evaluate the limit $s\rightarrow 1$. The two summation terms in \eqref{Cor4} follow immediately with the help of \eqref{Cor1} and the arguments used in its derivation. It remains to show that\newline

{\bf Lemma}
% from corollary4a.mw

\begin{equation}
\displaystyle  \lim _{s\rightarrow 1}\left( 2\,\ln  \left( 2 \right) S \left( s,q \right) -S^{\prime}   \left( s,q \right)\right) =2\sum _{j=0}^{\infty }{\frac { \left( -1 \right) ^{j}\ln  \left( j+2\,q \right) }{ \left( j+2\,q \right)}}
\mbox{}+\ln  \left( 2 \right)  \left( \psi \left( q+1/2 \right) -\psi \left( q \right)  \right) 
\label{Zids}
\end{equation}

\iffalse
By definition,  the Stieltjes constants $\gamma_{n}(q)$ are the coefficients of the Laurent expansion of $\zeta(s,q)$ about $s=1$, that is
\begin{equation}
\displaystyle \zeta \left(s,q \right) = \left( s-1 \right) ^{-1}+\sum _{n=0}^{\infty }{\frac { \left( -1 \right) ^{n}\gamma_{n} \left( a \right)  \left( s-1 \right) ^{n}}{\Gamma \left( n+1 \right) }}\,.
\label{Stieltje1}
\end{equation}
After differentiating once with respect to $s$, we have
\begin{equation}
\displaystyle \zeta^{\prime} \left( s,q \right) =- \left( s-1 \right) ^{-2}+\sum _{n=0}^{\infty }{\frac { \left( -1 \right) ^{n+1}\gamma_{n+1} \left(q \right)  \left( s-1 \right) ^{n}}{\Gamma
\mbox{} \left( n+1 \right) }}
\label{StieltjesD}
\end{equation}

We are interested in 
\begin{equation}
\lim _{s\rightarrow 1}(\zeta^{\prime} \left( s,q \right)-\zeta^{\prime} \left( s,q+1/2 \right))=\gamma_{1}(q+1/2)-\gamma_1(q)\,.
\end{equation}
\fi
{\bf Proof:}

%Differentiate \eqref{SaDef} with respect to $s$, giving
From the first equality in \eqref{clarity} and the definition \eqref{ZaltDef} we have
\begin{equation}
\displaystyle \zeta^{\prime} \left(s,q \right) -\zeta^{\prime} \left(s,q+1/2 \right) ={2}^{s}\ln  \left( 2 \right) \sum _{j=0}^{\infty }{\frac { \left( -1 \right) ^{j}}{ \left( j+2\,q \right) ^{s}}}
\mbox{}-{2}^{s}\sum _{j=0}^{\infty }{\frac { \left( -1 \right) ^{j}\ln  \left( j+2\,q \right) }{ \left( j+2\,q \right) ^{s}}}.
\label{Zidd}
\end{equation}
Since the  first sum in \eqref{Zidd} is simply identified as
\begin{equation}
\displaystyle \sum _{j=0}^{\infty }{\frac { \left( -1 \right) ^{j}}{ \left( j+2\,q \right) ^{s}}}={2}^{-s}S \left( s,q \right) 
\label{DefSa}
\end{equation}
the requisite Lemma is proven by taking the limit $s\rightarrow1$ and taking note of \eqref{Lims1}. {\bf QED}\newline

\begin{itemize}
\item

{\bf Special case n=0.}
% from file Cor4Special cases.mw

\begin{align} \nonumber
\displaystyle a\int_{0}^{\infty }\!\frac {v\arctan \left( {\frac {2v}{a}} \right)}{ \left( {a}^{2}+4\,{v}^{2} \right) \cosh \left( \pi\,v \right) } \,{\rm d}v&-\int_{0}^{\infty }\!{\frac {{v}^{2}\ln  \left( {a}^{2}+4\,{v}^{2} \right) }{ \left( {a}^{2}+4\,{v}^{2} \right) \cosh \left( \pi\,v \right) }}\,{\rm d}v\\ \nonumber
&=  a \ln  \left( 2 \right)\left( \psi \left( q+1/2 \right) -\psi \left( q \right)  \right)/8 +a/4\,\sum _{j=0}^{\infty }{\frac { \left( -1 \right) ^{j}\ln  \left( j+2\,q \right) }{j+2\,q}}
\mbox{}  \\
&+\left(\ln  \left( \Gamma \left( q \right)  \right) 
\mbox{}-\ln  \left( \Gamma \left( q+1/2 \right)  \right) -\ln  \left( 2 \right)\right)/2 
\label{C4an0}
\end{align}
\item
{\bf Special case n=0, a=1}

% from file Cor4Special cases.mw

The identity \cite{Math} 
\begin{equation}
\displaystyle \sum _{j=0}^{\infty }{\frac { \left( -1 \right) ^{j}\ln  \left( j+1 \right) }{j+1}}=\ln  \left( 2 \right)  \left( \ln  \left( 2 \right) -2\,\gamma \right)/2 
\label{Zaltm1}
\end{equation}
leads to
\begin{equation}
\displaystyle \int_{0}^{\infty }\!{\frac {v\arctan \left( 2\,v \right) }{ \left( 4\,{v}^{2}+1 \right) \cosh \left( \pi\,v \right) }}\,{\rm d}v-\int_{0}^{\infty }\!{\frac {{v}^{2}\ln  \left( 4\,{v}^{2}+1 \right) }{ \left( 4\,{v}^{2}+1 \right) \cosh \left( \pi\,v \right) }}\,{\rm d}v
=\frac{1}{4}\left(\frac{3}{2}  \ln^{2}  \left( 2 \right) -\gamma\ln  \left( 2 \right) +\ln  \left( \pi/4 \right)\right); 
\label{C4n0a1}
\end{equation}
\item
{\bf Special case n=0, a=2}

With recourse to \eqref{C2na2} with $n=0$, we find

\begin{align} \nonumber
\displaystyle 2\,\int_{0}^{\infty }\!{\frac {v\arctan \left( v \right) }{ \left( {v}^{2}+1 \right) \cosh \left( \pi\,v \right) }}\,&{\rm d}v-\int_{0}^{\infty }\!{\frac {{v}^{2}\ln  \left( {v}^{2}+1 \right) }{ \left( {v}^{2}+1 \right) \cosh \left( \pi\,v \right) }}\,{\rm d}v\\
&=4\,\sum _{j=0}^{\infty }{\frac { \left( -1 \right) ^{j}\ln  \left( 2\,j+3 \right) }{2\,j+3}}+\pi\,\ln  \left( 2 \right) +4\,\ln  \left( \Gamma \left( 3/4 \right)  \right) 
\mbox{}-2\,\ln  \left( 2\,\pi \right) \,;
\label{C4n0a2}
\end{align}

\item
{\bf Special case n=0, a=4}

With recourse to \eqref{C2na4} with $n=0$, we find
\begin{align} \nonumber
&\displaystyle 4\,\int_{0}^{\infty }\!{\frac {v\arctan \left( v/2 \right) }{ \left( {v}^{2}+4 \right) \cosh \left( \pi\,v \right) }}\,{\rm d}v-\int_{0}^{\infty }\!{\frac {{v}^{2}\ln  \left( {v}^{2}+4 \right) }{ \left( {v}^{2}+4 \right) \cosh \left( \pi\,v \right) }}\,{\rm d}v
\\
&=8\,\sum _{m=0}^{\infty }{\frac { \left( -1 \right) ^{m}\ln  \left( 2\,m+5 \right) }{2\,m+5}}-2\,\pi\,\ln  \left( 2 \right) -4\,\ln  \left( \Gamma \left( 3/4 \right)  \right) 
\mbox{}+16\ln  \left( 2 \right)/3 +2\,\ln  \left( \pi/3 \right)\,. 
\label{C4n0a4}
\end{align}

({\bf Remark}: Further evaluations are accessible by differentiating any of the above with respect to the parameter $a$.)
\end{itemize}

\section{Theorem 2 and Proof} \label{sec:T2andProof}

{\bf Theorem 2}
% from file "...\2016\new_master_fumctions\redo_of_old_paperfrom 2014\check_new_paper.mw"

\begin{equation}
\displaystyle \int_{0}^{\infty }\!{\frac {{v}^{2\,n}\,\cos \left( s\arctan \left( {\frac {2v}{a}} \right)  \right)}{ \left( {a}^{2}+4\,{v}^{2} \right) ^{s/2}\cosh \left( \pi\,v \right) 
} }\,{\rm d}v={2}^{-s} \left( -1 \right) ^{n} \left( \frac{a}{2} \right) ^{2\,n}\Gamma \left( 2\,n+1 \right) \sum _{j=0}^{2\,n} \frac {\left(-2/a \right) ^{j}\eta \left( s-j,2\,q \right) }{\Gamma \left( 2\,n-j+1 \right) \Gamma \left( j+1 \right) }
\label{Theorem2}
\end{equation}
\iffalse
\begin{equation}
\displaystyle \int_{0}^{\infty }\!{\frac {{v}^{2\,n}\cos \left( s\arctan \left( {\frac {2\,v}{a}} \right)  \right)  }{\left( {a}^{2}+4\,{v}^{2} \right) ^{s/2}\cosh \left( \pi\,
\mbox{}v \right) }}\,{\rm d}v={2}^{2(n-s)} \left( -1 \right) ^{n}\Gamma \left( 2\,n+1 \right) \sum _{j=0}^{2\,n}{\frac { \left( -a/4 \right) ^{j}S \left( j+s-2\,n,q \right) }{\Gamma \left( 2\,n-j+1 \right) \Gamma \left( j+1 \right) }}
\label{Theorem2}
\end{equation}
\fi
% for all of this See file Patk.mw and better yet, Theorem2.mw and Theorem2a.mw

{\bf Proof:}\newline

From \cite[Eq. 4.111(7)]{G&R}, after invocation of minor arctangent identities, we have the Fourier sine transform 
\begin{equation}
\displaystyle \int_{0}^{\infty }\!{\frac {v^{-1}\sin \left( wv \right) }{\cosh \left( \pi\,v/2 \right)}}\,{\rm d}v=2\,\arctan \left( \tanh
\mbox{} \left( w/2 \right)  \right)\,. 
\label{FTnm1a}
\end{equation}
In complete analogy to \eqref{Thm4a} this leads to

\begin{equation}
\displaystyle \int_{0}^{\infty }\!{\frac {{v}^{2\,n-1}\sin \left( s\arctan \left( {\frac {v}{a}} \right)  \right) }{ \left( {a}^{2}+{v}^{2} \right) ^{s/2}\cosh \left(\pi\,v/2 \right) 
\mbox{}}}\,{\rm d}v=2\, \left( -1 \right) ^{n}\Gamma \left( 2\,n+1 \right) \sum _{j=0}^{2\,n}{\frac { \left( -a \right) ^{j}J \left( j+s-2\,n-1 \right) }{\Gamma \left( 2\,n-j+1 \right) \Gamma \left( j+s-2\,n \right) \Gamma \left( j+1 \right) }}\,,
\label{Jid}
\end{equation}
where
\begin{equation}
\displaystyle J \left( s \right)\equiv \int_{0}^{\infty }\!{w}^{s}{{\rm e}^{-aw}}\arctan \left( \tanh \left( w/2 \right)  \right) \,{\rm d}w
\label{J(s)}
\end{equation}
obeys the recursion
\begin{equation}
\displaystyle J \left( s \right) ={\frac {sJ \left(s -1 \right) }{a}}+{\frac {{2}^{-1-2\,s}\Gamma \left( 1+s \right) S \left( 1+s,q \right) }{2\,a}}
\label{JgRecur}
\end{equation}
derived by integration by parts and the use of \eqref{Eq3p3}. Substitution of \eqref{JgRecur} into \eqref{Jid} with $s\rightarrow j+s-2n-1$ leads to

\iffalse
\begin{align} \nonumber
&\displaystyle \int_{0}^{\infty }\!{\frac {{v}^{2\,n-1}\sin \left( s\arctan \left( {\frac {v}{a}} \right)  \right)}{ \left( {a}^{2}+{v}^{2} \right) ^{s/2}\cosh \left( \pi\,v/2 \right) 
} }\,{\rm d}v =\frac { 2\,\left( -1 \right) ^{n}\Gamma \left( 2\,n+1 \right) }{a} \\
&\times\left( \sum _{j=0}^{2\,n}{\frac { \left( -a \right) ^{j}{\it J} \left( -2+j+s-2\,n \right) }{\Gamma \left( 2\,n-j+1 \right) \Gamma \left( j+s-2\,n-1 \right) 
\mbox{}\Gamma \left( j+1 \right) }} \right.   
 \left. +{2}^{-2\,s+4\,n}\sum _{j=0}^{2\,n}{\frac { \left( -a/4 \right) ^{j}S \left( j+s-2\,n,q \right) }{\Gamma \left( 2\,n-j+1 \right) \Gamma \left( j+1 \right) }}
\mbox{} \right) \,.
\label{T4a}
\end{align}
\fi

\begin{align} \nonumber
&\displaystyle \int_{0}^{\infty }\!{\frac {{v}^{2\,n-1}\sin \left( s\arctan \left( {\frac {v}{a}} \right)  \right)}{ \left( {a}^{2}+{v}^{2} \right) ^{s/2}\cosh \left( \pi\,v/2 \right) 
} }\,{\rm d}v =\frac { 2\,\left( -1 \right) ^{n}\Gamma \left( 2\,n+1 \right) }{a} \\
&\times\left( \sum _{j=0}^{2\,n}{\frac { \left( -a \right) ^{j}{\it J} \left( -2+j+s-2\,n \right) }{\Gamma \left( 2\,n-j+1 \right) \Gamma \left( j+s-2\,n-1 \right) 
\mbox{}\Gamma \left( j+1 \right) }}  +\, {2}^{2\,n-s} \left( \frac{a}{2} \right) ^{2\,n}\sum _{j=0}^{2\,n}{\frac {\left( -2/a \right) ^{j}\eta \left( s-j,2\,q \right) }{\Gamma \left( 2\,n-j+1 \right) \Gamma \left( j+1 \right) } }\right).
\label{T4a}
\end{align}

Now let $s:=s-1$ in \eqref{Jid} to find
\begin{align} \nonumber
\displaystyle \int_{0}^{\infty }\!\frac {{v}^{2\,n-1}\sin \left(  \left( s-1 \right) \arctan \left( {\frac {v}{a}} \right)  \right)}{ \left( {a}^{2}+{v}^{2} \right) ^{-1/2+s/2}
\mbox{}\cosh \left(\pi\,v/2 \right) }& \,{\rm d}v\\
&\hspace{-80pt}=2\, \left( -1 \right) ^{n}\Gamma \left( 2\,n+1 \right) \sum _{j=0}^{2\,n}{\frac { \left( -a \right) ^{j}J \left( -2+j+s-2\,n \right) }{\Gamma \left( 2\,n-j+1 \right) \Gamma \left( j+s-2\,n-1 \right) \Gamma \left( j+1 \right) }}
\label{T4A}
\end{align}

Notice that the sum on the right-hand side of \eqref{T4A} coincides with the first sum on the right-hand side of \eqref{T4a}, so substitute, expand the term ${ { \sin(\left(  s-1 \right) \arctan \left( { {v}/{a}} \right)) }}$ using a simple trigonometric identity and simplify. The result is \eqref{Theorem2}. {\bf QED}\newline

The above generalizes the classic result (see \cite[Eq. (26)]{Milgram} and references therein) corresponding to $n=0,a=1$ and Glasser's more recent generalization \cite{Glasserlattice1}, corresponding to the case $n=0$. Additionally, various specific cases extend results found in \cite[Section 3.522]{G&R}.\newline

{\bf Special Cases}
\begin{itemize}

\item
Case: $s=0$.
When reduced by use of \eqref{Zeta2Bern} and \eqref{BernId1}, this case simply reproduces a known result \cite[Eq. 3.523(4)]{G&R}.

\item

Case: $s=1$.

After an obvious change of variables using $n=0$ and $a=1/2$ we find
\begin{equation}
\displaystyle \int_{0}^{\infty }\!{\frac {1}{ \left( {v}^{2}+1 \right) \cosh \left( \pi\,v/4 \right) }}\,{\rm d}v=\left(\psi \left( 7/8 \right) -\psi \left( {{3}/{8}} \right) \right)/2
\label{T1s0}
\end{equation}
\end{itemize}
Comparison with \cite[Eq. 3.522(10)]{G&R} yields the identity
\begin{equation}
\displaystyle \psi \left( 7/8 \right) -\psi \left( {{3}/{8}} \right) = \sqrt{2} \left( \pi-2\,\ln  \left(  \sqrt{2}+1 \right)  \right)\,,	 
\label{PsiId}
\end{equation}
which could alternatively have been obtained from Gauss' classic expression (see \cite[Eq. (1.3)]{Choi&Cvij}).

\subsection{Corollary 4.1} \label{sec:4.1}

% This is from file Theorem3 where it is derived in a far more complicated manner by adding and subtracting 4.2 and 3.4

\begin{align}
\displaystyle \int_{0}^{\infty }\!{\frac {{v}^{2\,n}\ln  \left( {a}^{2}+4\,{v}^{2} \right) }{\cosh \left( \pi\,v \right) }}\,{\rm d}v&= \left( -1 \right) ^{n}  \biggl(  2 ^{1-2\,n}\mathcal{E} \left( 2\,n \right) \ln  \left( 2 \right)- {2}^{2\,n+1}\,S^{\prime}(-2n,q) \\
&  -2\, \left( a/2 \right) ^{2\,n}\Gamma \left( 2\,n+1 \right) \sum _{j=0}^{2\,n-1}{\frac {\left( -4/a \right) ^{j} S^{\prime}(-j,q) }{\Gamma \left( 2\,n-j+1 \right) \Gamma \left( j+1 \right) 
\mbox{}} }
\mbox{}  \biggr)
\label{Lint5}
\end{align}
{\bf Proof:} Operate on \eqref{Theorem2} with $\frac{\partial}{\partial s}$ and let $s=0$. The proof follows exactly the same steps as that of Section \ref{sec:3.4}.   {\bf QED}
\begin{itemize}
\item
Case $n=0$: (also see \cite[Eqs. 4.373(1) and (2)]{G&R})

\begin{equation}
\displaystyle \int_{0}^{\infty }\!{\frac {\ln  \left( {a}^{2}+4\,{v}^{2} \right) }{\cosh \left( \pi\,v \right) }}\,{\rm d}v=2\,\ln  \left( 2 \right) -2\,\ln  \left( \Gamma \left( q \right)  \right) 
\mbox{}+2\,\ln  \left( \Gamma \left( q+1/2 \right)  \right) 
\label{C4n0}
\end{equation}

\item
Case: $n=1$

% from file Theorem3.mw

\begin{equation}
\displaystyle \int_{0}^{\infty }\!{\frac {{v}^{2}\ln  \left( {a}^{2}+4\,{v}^{2} \right) }{\cosh \left( \pi\,v \right) }}{\rm d}v=\frac{a^{2}}{2}\left(\ln  \left( \Gamma \left( q \right)  \right) -\ln  \left( \Gamma \left( q+1/2 \right)  \right)   \right) 
\mbox{}-4a\,S^{\prime} \left( -1,q \right) +8\,S^{\prime} \left( -2,q \right) 
\mbox{}+\frac{\ln  \left( 2 \right)}{2} 
\label{C4n1}
\end{equation}
\item
Case $n=1, a=1$\newline

From \eqref{Z1m1}, \eqref{Z1p2}, \eqref{Zm1h}, \eqref{Zm2h} and\eqref{Zpm2}
\begin{equation}
\displaystyle \int_{0}^{\infty }\!{\frac {{v}^{2}\ln  \left( 4\,{v}^{2}+1 \right) }{\cosh \left( \pi\,v \right) }}\,{\rm d}v=\ln  \left( \pi \right)/4 +2\,\ln  \left( 2 \right)/3 +1/2
\mbox{}-6\,\ln  \left( A \right) +{\frac {7\,\zeta \left( 3 \right) }{{2\,\pi}^{2}}}
\label{C4n1a0}
\end{equation}

\item
% from file Hyper_recursion.mw
Case $n=2, a=1$ \newline
\begin{equation}
\displaystyle \int_{0}^{\infty }\!{\frac {{v}^{4}\ln  \left( 4\,{v}^{2}+1 \right) }{\cosh \left( \pi\,v \right) }}\,{\rm d}v={\frac {23\,\ln  \left( 2 \right) }{40}}+{\frac {93\,\zeta \left( 5 \right) }{2\,{\pi}^{4}}}
\mbox{}-60\,\zeta^{\prime} \left(-3 \right) -{\frac {21\,\zeta \left( 3 \right) }{4\,{\pi}^{2}}}-\frac{1}{4}+3\,\ln  \left( A \right) -\frac{\ln(\pi)}{16}. 
\label{R1}
\end{equation}
\end{itemize}
\subsection{Corollary 4.2} \label{sec:4.2}

% see file Theorem2a.mw

\begin{align} \nonumber
&\displaystyle \int_{0}^{\infty }\!{\frac {{v}^{2\,n+1}\arctan \left( {\frac {2v}{a}} \right)}{ \left( {a}^{2}+4\,{v}^{2} \right) \cosh \left( \pi\,v \right) } }\,{\rm d}v+\frac{a}{4}\int_{0}^{\infty }\!{\frac {{v}^{2\,n}\ln  \left( {a}^{2}+4\,{v}^{2} \right) }{ \left( {a}^{2}+4\,{v}^{2} \right) \cosh \left( \pi\,v \right) }}\,{\rm d}v
\mbox{}={2}^{2\,n-3} \left( -1 \right) ^{n} \left( a/4 \right) ^{2\,n}\\ \nonumber
& \left( -{\frac {4\,\Gamma \left( 2\,n
\mbox{}+1 \right) }{a} \left( \ln  \left( 2 \right) \sum _{j=0}^{2\,n-1}{\frac {\mathcal{E}\left( j,2\,q \right)\left( -2/a \right) ^{j} }{\Gamma \left( j+2 \right) \Gamma \left( 2\,n-j \right) } }
\mbox{}-\sum _{j=0}^{2\,n-1}{\frac { \left( -4/a \right) ^{j}\,S^{\prime}  \left( -j\,,q \right) }{\Gamma \left( j+2 \right) \Gamma \left( 2\,n-j \right) }} \right) } \right. \\
&\left.\hspace{100pt}+2\,\sum _{j=0}^{\infty }{\frac { \left( -1 \right) ^{j}\ln  \left( j+2\,q \right) }{j+2\,q}}+\ln  \left( 2 \right)  \left( \psi \left( q+1/2 \right) -\psi \left( q \right)  \right) 
\mbox{} \right) \,.
\label{T2s1D}
\end{align}

{\bf Proof:} Operate on \eqref{Theorem2} with $\frac{\partial}{\partial s}$ and let $s\rightarrow 1$. The proof follows exactly the same steps as that of Section \ref{sec:3.4}.   {\bf QED}

\section{A Comment on the Even moments of Theorem 1} \label{sec:Even}

As noted previously, in his paper, Patkowski purports to resolve the case where the identity \eqref{Theorem4} contains even rather than odd moments. As discussed in the Appendix, the derivation of that result is flawed; additionally the putative result does not satisfy numerical testing. To utilize the method developed by Patkowski, the Fourier sine transform of ${ t^{2n}/\cosh(at)}$ for {\bf any} $n\geq 0$, must be employed, and in fact, such an identity is known for $n=0$. From \cite[Eq. 3.981(2)]{G&R} we find
% from file "Even_moments.mw"

\begin{equation}
\displaystyle \int_{0}^{\infty }\!{\frac {\sin \left( wx \right) }{\cosh \left(\pi\,x/2 \right) }}\,{\rm d}x={\frac {i}{\pi}\psi \left( \frac{1}{4}-{\frac {iw}{2\pi}} \right) }
\mbox{}-{\frac {i}{\pi}\psi \left( {\frac{1}{4}+\frac {iw}{2\pi}} \right) }-\tanh \left( w \right) ,
\label{G&Reven}
\end{equation}
which may also be written as

\begin{equation}
\displaystyle \int_{0}^{\infty }\!{\frac {\sin \left( wx \right) }{\cosh \left( \pi\,x/2 \right) }}\,{\rm d}x=-\tanh \left( w \right) +{\frac {2}{\pi}\Im \left( \psi \left( {\frac{1}{4}+\frac {iw}{2\pi}} \right)  \right) }\,.
\label{E3}
\end{equation}
From Parseval's identity, we then have

\begin{align}
\displaystyle \int_{0}^{\infty }&\!{\frac {{t}^{2\,n}\sin \left( s\arctan \left( {\frac {t}{a}} \right)  \right)}{ \left( {a}^{2}+{t}^{2} \right) ^{s/2}\cosh \left( \pi\,t/2 \right) 
\mbox{}} }\,{\rm d}t=-\frac { \left( -1 \right) ^{n}\Gamma \left( 2\,n+1 \right) }{\pi\,\Gamma \left( s \right) }\left(\pi\,J_1(n,s,a)-2\,J_2(n,s,a)\right)
\label{E4a}
\end{align}
where
\begin{equation}
J_{1}(n,s,a)\equiv \int_{0}^{\infty }\!{L} \left( 2\,n,s-2\,n-1,aw \right) {w}^{s-2\,n-1}
\mbox{}\tanh \left( w \right) {{\rm e}^{-aw}}\,{\rm d}w
\label{J1}
\end{equation}
and
\begin{equation}
J_{2}(n,s,a)\equiv \int_{0}^{\infty }\!{L} \left( 2\,n,s-2\,n-1,aw \right) {w}^{s-2\,n-1}\Im \left( \psi \left( \frac{1}{4}\,+\frac{iw}{2\pi} \right)  \right)
\mbox{}{{\rm e}^{-aw}}\,{\rm d}w\,.
\label{J2}
\end{equation}
From \eqref{LagDef} these become
\begin{equation}
\displaystyle J_{1}(n,s,a)=\Gamma \left( s \right) \sum _{j=0}^{2\,n} \left( {\frac { (-a)^{j}\int_{0}^{\infty }\!{w}^{j+s-2\,n-1}\tanh \left( w \right) {{\rm e}^{-aw}}\,{\rm d}w}{\Gamma \left( 2\,n-j+1 \right) 
\mbox{}\Gamma \left( j+s-2\,n \right) \Gamma \left( j+1 \right) }} \right) 
\label{J1a}
\end{equation}
and
\begin{equation}
\displaystyle J_{{2}} \left( n,s,a \right) =\Gamma \left( s \right) \sum _{j=0}^{2\,n}  {\frac { (-a)^{j}\int_{0}^{\infty }\!{w}^{j+s-2\,n-1}\Im \left( \psi \left( {\frac{1}{4}+\frac {iw}{2\pi}} \right)  \right) 
\mbox{}{{\rm e}^{-aw}}\,{\rm d}w}{\Gamma \left( 2\,n-j+1 \right) \Gamma \left( j+s-2\,n \right) \Gamma \left( j+1 \right) }} \,.
\label{J2a}
\end{equation}
With $\Re(s)>-2$, the first of these leads to the related integral
\begin{align} \nonumber
\displaystyle \int_{0}^{\infty }\!{w}^{s}\tanh \left( w \right) {{\rm e}^{-aw}}\,{\rm d}w&={2}^{-2-2\,s}\Gamma \left( 1+s \right) 
\mbox{} \\
&\times\left( \zeta \left( 1+s,1+a/4 \right) +\zeta \left(1+s,a/4 \right) -2\,\zeta \left(1+s,1/2+a/4 \right)  \right)\,, 
\label{Jint1}
\end{align}

the equality arising courtesy of Mathematica \cite{Math}, verified numerically. The evaluation of the integral in \eqref{J2a} is more challenging. With $\Re(s)>-2$, we require an evaluation of the generic form
\begin{equation}
J_{2}(s,a)\equiv\int_{0}^{\infty }\!{w}^{s}\Im \left( \psi \left( {\frac{1}{4}+\frac {iw}{2\pi}} \right)  \right) 
\mbox{}{{\rm e}^{-aw}}\,{\rm d}w\,.
\label{Jint2}
\end{equation}

One possibility is to identify
% from file PsId.mw

\begin{equation}
\displaystyle \Im \left( \psi \left( {\frac{1}{4}+\frac{iw}{2\pi}} \right)  \right) ={\frac {1}{2\pi}\sum _{k=1}^{\infty }\frac{w}{  \left( k-3/4 \right) ^{2}
\mbox{}+{\displaystyle  {{w}^{2}}/{{4\pi}^{2}}}  }}\,.
\label{PsId}
\end{equation}
Substituted into \eqref{Jint2}, this leads to known identities (see \cite[Eqs. 3.356(1) and 3.356(2)] {G&R}), but only when $s$ is a positive integer. The case of arbitrary $s$ does not appear to be tractable. Alternatively, and following integration by parts, \eqref{Jint2} leads to an interest in the real-valued integral

\begin{align}
\displaystyle \overset{-}{J}(s,a)&\equiv\,\int_{0}^{\infty }\!{w}^{s}{{\rm e}^{-aw}} \left( \log  \left( \Gamma \left( \frac{1}{4}+{\frac {iw}{2\pi}} \right)  \right) +\log  \left( \Gamma \left( \frac{1}{4}-{\frac {iw}{2\pi}} \right)  \right) 
\mbox{} \right) \,{\rm d}w\\
&=2\,\int_{0}^{\infty }\!{w}^{s}{{\rm e}^{-aw}} \log  \left| \Gamma \left( \frac{1}{4}+{\frac {iw}{2\pi}} \right)  \right| 
\mbox{}  \,{\rm d}w
\label{Jint2b}
\end{align}
 
because
\begin{equation}
J_{2}(s,a)=\pi s \overset{-}{J}(s-1,a)-a\overset{-}{J}(s,a)\,.
\end{equation}
And this integral appears also to be intractable. Therefore, further analysis of the case under consideration in this Section must be left unresolved.

\section{Connection with $\zeta(2m+1)$} \label{sec:Eta}

\iffalse
Noting that, for $s=2m+1$, Dirichlet's function $\eta(s)$ is related to Riemann's function $\zeta(s)$ by
\begin{equation}
\eta(2m+1)= \left( {2}^{2\,m+1}-1 \right) \zeta \left( 2\,m+1 \right) \,,
\end{equation}
\fi

From \eqref{Cp5a} and \eqref{A(m)}, we are interested in

\begin{equation}
\displaystyle A \left( m \right) ={\frac {{2}^{3+2\,m}
\mbox{}{\pi}^{2\,m}}{\Gamma \left( 3+2\,m \right) }\int_{0}^{\infty }\!{\frac {{t}^{2\,m+2}{\mbox{$_3$F$_2$}(1,1,3/2;\,2+m,3/2+m;\,-4\,{t}^{2})}}{\cosh \left( \pi\,t \right) }}\,{\rm d}t}\,.
\end{equation} 

% see file ..\2019\EiSum\hyper_recursion.mw

From \cite[Eq. (3)] {Shpot&Sriv}
\iffalse
where
\begin{align} \nonumber
&\displaystyle {\mbox{$_3$F$_2$}(a,b,c;\,b+1+m,c+1+n;\,z)}=\frac {\Gamma \left( b+1+m \right) \Gamma \left( c+1+n \right) }{  \Gamma \left( b \right) \Gamma \left( c \right) } \\ \nonumber
&\times \sum _{k=0}^{m} \sum _{j=0}^{n}\frac { \left( -1 \right) ^{k}
\mbox{}  \left( -1 \right) ^{j} }{\Gamma \left( 1-k+m \right) \Gamma \left( n+1-j \right) \Gamma \left( k+1 \right) 
\mbox{}\Gamma \left( j+1 \right) }  \\  \label{Shpot1}
&\hspace{3cm}\times \left( {\frac {{\mbox{$_2$F$_1$}(a,c+j;\,c+j+1;\,z)}}{ \left( c+j \right)  \left( b-c+k-j \right) }}+{\frac {{\mbox{$_2$F$ _1$}(a,b+k;\,b+k+1;\,z)}}{ \left( b+k \right)  \left( c-b+j-k \right) }} \right)   
\end{align}
\fi
we obtain

\begin{align} \nonumber
&\displaystyle {\mbox{$_3$F$_2$}(1,1,3/2;\,2+m,m+3/2;\,z)}=\frac {2\,\Gamma \left( 2+m \right) \Gamma \left( m+3/2 \right) }{ \sqrt{\pi}}\\ \nonumber
&\times\sum _{k=0}^{m}  \sum _{j=0}^{m-1}\frac { \left( -1 \right) ^{k} \left( -1 \right) ^{j}}{\Gamma \left( 1-k+m \right) 
\mbox{}\Gamma \left( -j+m \right) \Gamma \left( k+1 \right) \Gamma \left( j+1 \right) }\\
&\hspace{3cm}\times \left( {\frac {{\mbox{$_2$F$_1$}(1,3/2+j;\,5/2+j;\,z)}}{ \left( 3/2+j \right)  \left( -1/2+k-j \right) }}
\mbox{}+{\frac {{\mbox{$_2$F$_1$}(1,k+1;\,2+k;\,z)}}{ \left( k+1 \right)  \left( 1/2+j-k \right) }} \right) \,.  
\label{Shpot2}
\end{align}
Further, from \cite[Eq. 7.3.1(134) and (135)]{prudnikov} with $t\geq 0$ we respectively have
\begin{equation}
\displaystyle {\mbox{$_2$F$_1$}(1,3/2+j;\,5/2+j;\,-4\,{t}^{2})}={\frac { \left( -1 \right) ^{j+1}
\mbox{} \left( 2\,j+3 \right) }{ \left( 4\,{t}^{2} \right) ^{j+2}} \left( 2\,t\arctan \left( 2\,t \right) +\sum _{p=1}^{j+1}{\frac { \left( -4\,{t}^{2} \right) ^{p}}{2\,p-1}} \right) }
\label{H1}
\end{equation}
and
\begin{equation}
\displaystyle {\mbox{$_2$F$_1$}(1,k+1;\,2+k;\,-4\,{t}^{2})}=-{\frac {k+1}{ \left( -4\,{t}^{2} \right) ^{k+1}} \left( \ln  \left( 4\,{t}^{2}+1 \right) +\sum _{p=1}^{k}{\frac { \left( -4\,{t}^{2} \right) ^{p}}{p}} \right) }\,.
\label{H2}
\end{equation}
For the case $m=1$, \eqref{Cp5a} and the above lead to the representation

% again from file hyper_recursion.mw
\begin{align} \nonumber
\displaystyle {\frac {\zeta \left( 3 \right) }{{\pi}^{2}}}=\frac{2}{7}\,\int_{0}^{\infty }\!{\frac {{t}^{2}\ln  \left( 4\,{t}^{2}+1 \right) }{\cosh \left( \pi\,t \right) }}\,{\rm d}t
\mbox{}&-\frac{1}{14}\,\int_{0}^{\infty }\!{\frac {\ln  \left( 4\,{t}^{2}+1 \right) }{\cosh \left( \pi\,t \right) }}\,{\rm d}t+\frac{4}{7}\,\int_{0}^{\infty }\!{\frac {t\arctan \left( 2\,t \right) }{\cosh \left( \pi\,t \right) }}\,{\rm d}t
\mbox{}\\
& -\frac{6}{7}\,\int_{0}^{\infty }\!{\frac {{t}^{2}}{\cosh \left( \pi\,t \right) }}\,{\rm d}t+{\frac{3}{28}}\,.
\label{Casem1}
\end{align}
From \eqref{C4n1a0}, \eqref{C4n0}, \eqref{C3n0a1} and  \cite[Eq. 3.523(5)]{G&R} in that order, \eqref{Cp5a} is verified when $m=1$.\newline

Similarly, for the case $m=2$ we have
% see file hyper_recursion.mw
\begin{align} \nonumber
\displaystyle {\frac {\zeta \left( 5 \right) }{{\pi}^{4}}}&={\frac {2}{93}\int_{0}^{\infty }\!{\frac {\ln  \left( 4\,{t}^{2}+1 \right) {t}^{4}}{\cosh \left( \pi\,t \right) }}\,{\rm d}t}
\mbox{}-\frac{1}{31}\,\int_{0}^{\infty }\!{\frac {\ln  \left( 4\,{t}^{2}+1 \right) {t}^{2}}{\cosh \left( \pi\,t \right) }}\,{\rm d}t+{\frac {1}{744}\int_{0}^{\infty }\!{\frac {\ln  \left( 4\,{t}^{2}+1 \right) }{\cosh \left( \pi\,t \right) }}\,{\rm d}t}
\mbox{}\\ \nonumber
&+{\frac {8}{93}\int_{0}^{\infty }\!{\frac {\arctan \left( 2\,t \right) {t}^{3}}{\cosh \left( \pi\,t \right) }}\,{\rm d}t}-{\frac {2}{93}\int_{0}^{\infty }\!{\frac {t\arctan \left( 2\,t \right) }{\cosh \left( \pi\,t \right) }}\,{\rm d}t}
\mbox{}-{\frac {25}{279}\int_{0}^{\infty }\!{\frac {{t}^{4}}{\cosh \left( \pi\,t \right) }}\,{\rm d}t}\\
&+{\frac {7}{186}\int_{0}^{\infty }\!{\frac {{t}^{2}}{\cosh \left( \pi\,t \right) }}\,{\rm d}t}+{\frac{83}{8928}}\,.
\label{C2b}
\end{align}
All of the above integrals are given in the various examples presented in prior Sections, save for the second last, which is given in \cite[Eq. 3.523(7)]{G&R}. When the various substitutions are made, \eqref{Cp5a} is verified, this time with $m=2$, given that
\begin{equation}
\displaystyle \zeta^{\prime} \left(-4 \right) ={\frac {3\,\zeta \left( 5 \right) }{{4\,\pi}^{4}}}
\label{Zpminus4}
\end{equation}
and
\begin{equation}
\displaystyle \zeta^{\prime} \left(-4,1/2 \right) =-{\frac {45\,\zeta \left( 5 \right) }{64\,{\pi}^{4}}}\,.
\label{Zpm4half}
\end{equation}

\section{Summary}

A number of useful integral identities and sums have been obtained from a study of \eqref{Theorem4} and its special cases, utilizing Laguerre polynomials and Parseval's theorem. The important case discussed in Section \ref{sec:Even} remains unresolved.

\section{Acknowledgements}

The author thanks Larry Glasser for his helpful insights and independent derivation(s) of some of the special cases obtained here. All expenses related to this work have been borne by the author, who, not having an affiliation with a Canadian educational institution has been refused, by Canada's NSERC research council, permission to apply for a grant to allow Open Access publication. The ground for refusal is that the NSERC mandate only covers research grants for educational purposes at educational institutions. Therefore, if any student or staff at a Canadian educational institute happens to chance upon this paper, please do not read it.

\bibliographystyle{unsrt}

\bibliography{c://physics//biblio}

\begin{appendices}

\section{Patkowski's paper \cite{Patkowski2} Errata} \label{sec:AppendixA}
% all from file Patk.mw

The following are corrections to \cite{Patkowski2}; all numbered references to \cite{Patkowski2} in this Appendix are enclosed in square brackets (viz. ``[..]").\newline

Eq. [3.1] is missing a factor of two and should read

\begin{equation}
\displaystyle \int_{0}^{\infty }\!{\frac {{t}^{s-1}{{\rm e}^{-at}}}{\sinh \left( t \right) }}\,{\rm d}t=2\,\Gamma \left( s \right)  \left( \zeta \left(s,a \right) -{2}^{-s}\zeta \left(s,a/2 \right)  \right) 
\label{Eq3p1}
\end{equation}

Eq. [3.3] is missing a factor of two and should read as given in \eqref{Eq3p3}.\newline

Eq. [3.8] [i.e. proof of Theorem 2] has missing and extraneous factors and should read

\begin{align} \nonumber
\displaystyle \int_{0}^{\infty }\!&{\frac {{t}^{2\,n}\sin \left( s\arctan \left( {\frac {t}{a}} \right)  \right)}{ \left( {a}^{2}+{t}^{2} \right) ^{s/2}\sinh \left( \pi\,t \right) 
\mbox{}} }\,{\rm d}t\\
&={\frac { \left( -1 \right) ^{n}\Gamma \left( 2\,n+1 \right) \int_{0}^{\infty }\!{\it L} \left( 2\,n,s-2\,n-1,aw \right) \tanh \left( w/2 \right) 
\mbox{}{{\rm e}^{-aw}}{w}^{s-2\,n-1}\,{\rm d}w}{2\,\Gamma \left( s \right) }}\,.
\label{Eq3p8}
\end{align}

Eq. [3.4]  [i.e. Theorem 2] should read 
\begin{equation}
\displaystyle\int_{0}^{\infty }\!{\frac {{t}^{2\,n}\sin \left( s\arctan \left( {\frac {t}{a}} \right)  \right) }{ \left( {a}^{2}+{t}^{2} \right) ^{s/2}\sinh \left( \pi\,t \right) 
\mbox{}}}\,{\rm d}t=-\frac{a^{-s}}{2}\,\delta_{n,0}+\frac{1}{2}\,\sum _{m=0}^{2\,n} \left( -1 \right) ^{m+n}\binom{2\,n}{ m}{a}^{m}P_{{2}} \left( a,m+s-2\,n \right)\,,
\label{PatkThm2}
\end{equation}

where $\delta_{n,0}$ is the Kronecker delta, and $P_{2}(a,s)$ should have been defined as follows:

\begin{equation}
\displaystyle P_{{2}} \left( a,s \right) ={2}^{2-s}\zeta \left(s,a/2 \right) -2\,\zeta \left(s,a \right) \,.
\label{P2}
\end{equation}

Eq. [3.11] [ i.e. Theorem 3], should read
\begin{equation}
\displaystyle \int_{0}^{\infty }\!{\frac {{t}^{2\,n}\sin \left( s\arctan \left( {\frac {t}{a}} \right)  \right)}{ \left( {a}^{2}+{t}^{2} \right) ^{s/2} \left( {{\rm e}^{2\,\pi\,t}}+1 \right) 
\mbox{}} }\,{\rm d}t=\frac{1}{2}\sum _{m=0}^{2\,n} \left( -1 \right) ^{m+n}
\binom {2\,n}{ m}{a}^{m}P_{{3}} \left( a,m+s-2\,n \right) 
\label{PatTheorem3}
\end{equation}
where $P_{3}(a,s)$ should have been defined as follows:

\begin{equation}
\displaystyle P_{{3}} \left( a,s \right) ={\frac {{a}^{1-s}}{s-1}}+\zeta \left(s,a \right) -{2}^{s}\zeta \left(s,2\,a \right) \,.
\label{PatP3}
\end{equation}
{\bf Note:} The left hand side of \eqref{PatTheorem3} often appears in the literature written as

\begin{equation}
\displaystyle \int_{0}^{\infty }\!{\frac {{t}^{2\,n}\sin \left( s\arctan \left( {\frac {t}{a}} \right)  \right)}{ \left( {a}^{2}+{t}^{2} \right) ^{s/2} \left( {{\rm e}^{2\,\pi\,t}}+1 \right) 
\mbox{}} }\,{\rm d}t=\displaystyle \frac{1}{2}\,\int_{0}^{\infty }\!{\frac {{t}^{2\,n}\sin \left( s\arctan \left( {\frac {t}{a}} \right)  \right){{\rm e}^{-\pi\,t}}}{ \left( {a}^{2}+{t}^{2} \right) ^{s/2}
\mbox{}\cosh \left( \pi\,t \right) } }\,{\rm d}t\,.
\label{PatT3alt}
\end{equation}\newline

The derivation of Eq. [3.14]  [i.e. Theorem 4] employs the following correct identity
\begin{equation}
\displaystyle\int_{0}^{\infty }\!{\frac {\cos \left( wt \right) }{\cosh \left( \beta\,t \right) }}\,{\rm d}t=\frac{\pi}{2\,\beta\, \cosh \left( {\displaystyle \frac {\pi\,w}{2\,\beta}} \right)  }\,.
\label{PatkEq3p10}
\end{equation}
However, this identity is a Fourier {\bf cosine} transform, and thus cannot be used in Parseval's {\bf sine} identity as employed - see \eqref{Thm4a} - in \cite{Patkowski2}. A number of possible variations are listed in \cite[Eqs. 3.769 (1-7)]{G&R}; however none of these will yield a result different from those obtained here - see \eqref{Theorem4} or \eqref{Theorem2}. Thus Eq. [3.14] and Eq. [3.20] of \cite{Patkowski2}, based on \eqref{PatkEq3p10} are fundamentally flawed. See \eqref{Theorem4} for the correct result corresponding to [3.20] and Section \ref{sec:Even} for comments regarding Eq. [3.14] - i.e. [Theorem 4] of \cite{Patkowski2}. \newline

Eq. [3.17] would be correct if the right-hand side had the opposite sign. It should read
% see file "Check New paper.mw"
\begin{equation}
\displaystyle \int_{0}^{\infty }\!\frac {{t}^{2\,n+1}\cos \left( s\arctan \left( {\frac {t}{a}} \right)  \right)}{ \left( {a}^{2}+{t}^{2} \right) ^{s/2} \left( {{\rm e}^{2\,\pi\,t}}-1 \right) 
} \,{\rm d}t=-\frac{1}{2}\sum _{m=0}^{2\,n+1} \left( -1 \right) ^{m+n}\binom{2\,n+1}{ m}{a}^{m}{\it P}_{1} \left( a,m+s-2\,n-1 \right) \,,
\label{Pat3p17}
\end{equation}
where
\begin{equation}
\displaystyle {\it P}_{1} \left( a,s \right) =\zeta \left( s,a \right) -\frac{a^{-s}}{2}-{\frac {{a}^{1-s}}{s-1}}\,.
\label{PatP1}
\end{equation}
With the addition of a minus sign on the right-hand side, Eq. [3.18] of \cite{Patkowski2} would correctly read
\begin{equation}
\displaystyle \int_{0}^{\infty }\!{\frac {{t}^{2\,n+1}\cos \left( s\arctan \left( {\frac {t}{a}} \right)  \right)}{ \left( {a}^{2}+{t}^{2} \right) ^{s/2}\sinh \left( \pi\,t \right) 
\mbox{}} }\,{\rm d}t=-\frac{1}{2}\,\sum _{m=0}^{2\,n+1} \left( -1 \right) ^{m+n}\binom{2\,n+1}{m}{a}^{m}{\it P}_{2} \left( a,m+s-2\,n-1 \right) 
\label{Patk3p18}
\end{equation}
where $P_2$ is defined in \eqref{P2}.
%__________________________________
\end{appendices}
\end{flushleft}
%\end{maplegroup}
\end{document}